\documentclass[11pt]{article}
\usepackage{amssymb,amsmath,amsthm,epsfig}
\usepackage{latexsym, enumerate}
\usepackage{eepic}
\usepackage{epic}
\usepackage{graphicx}
\usepackage{color}
\usepackage{ifpdf}
\usepackage{subfigure}
\usepackage{tikz}
\usepackage{dsfont}
\usepackage{multirow}
\usepackage{makecell}
\usepackage{algorithm}
\usepackage{bm}
\usepackage{multirow}
\usepackage{epstopdf}

\usepackage{cite}
\usepackage{hyperref}
\usepackage{xcolor}
\hypersetup{
  colorlinks,
  linkcolor={blue!80!black},
  urlcolor={blue!80!black},
  citecolor={blue!80!black}
}
\usepackage[capitalize,noabbrev]{cleveref}
\usepackage{amsmath}
\numberwithin{equation}{section}

\theoremstyle{plain}
\newtheorem{lem}{Lemma}[section]
\newtheorem{thm}[lem]{Theorem}

\theoremstyle{definition}
\newtheorem{defn}{Definition}[section]

\theoremstyle{remark}
\newtheorem{rem}{Remark}[section]

\numberwithin{equation}{section}

\begin{document}
\title{
Determining internal topological structures and running cost of mean field games with partial boundary measurement}

\author{
Ming-Hui Ding\thanks{School of Mathematics and Statistics, Northwestern Polytechnical University, Xi'an, Shaanxi Province, China.\ \  Email: dingmh@nwpu.edu.cn}
\and
Hongyu Liu\thanks{Department of Mathematics, City University of Hong Kong, Kowloon, Hong Kong, China.\ \ Email: hongyu.liuip@gmail.com; hongyliu@cityu.edu.hk}
\and
Guang-Hui Zheng\thanks{School of Mathematics, Hunan University, Changsha 410082, China.\ \ Email: zhenggh2012@hnu.edu.cn; zhgh1980@163.com}
}

\date{}
\maketitle
\begin{center}{\bf ABSTRACT}
\end{center}\smallskip
The simultaneous recovery of the interior topological boundary and the nonlinear running cost functional from limited boundary observations remains a notable challenge in the inverse theory of mean-field game (MFG) systems. The problem presents substantial analytical difficulties due to the strongly coupled forward-backward parabolic structure of MFG models, the measure constraints satisfied by the population density, and the practical limitation that only partial boundary data are typically available in real-world applications. This work tackles this problem through a unified higher-order linearization framework tailored for MFG inverse problems under both Dirichlet and Neumann external boundary measurements. Rigorous uniqueness results are first established for the simultaneous identification of the cost functional and the interior geometric configuration under homogeneous boundary conditions. Building on this analytical foundation, we further extend the uniqueness theory to the practically more meaningful setting of inhomogeneous boundary data, substantially broadening the scope of applicable scenarios. The findings contribute significantly to the understanding of simultaneous reconstruction in complex MFG systems.

\smallskip


\section{Introduction}
\subsection{Background and Problem setup}
The study of mean-field game inverse problems has emerged as a highly anticipated research direction within the field of game theory. Since the groundbreaking work by Caines, Huang, Malhame, Lasry, and Lions in 2006-2007 \cite{Lasry2006a,Huang2006,Lasry2006,Huang2007,Lasry2007}, mean-field game theory has provided valuable tools and methodologies for understanding the equilibrium and behavior of large-scale intelligent systems. By assuming that each individual in a large-scale intelligent system focuses solely on the overall average effect when interacting with others, mean-field games simplify the analysis of complex systems. Research on inverse problems arising from mean-field dynamical systems possesses important theoretical significance and practical application value, and has attracted extensive attention across multiple disciplines. In economic research, for example, exploring inverse issues concerning market competition mechanisms and resource distribution strategies facilitates improved interpretation and forecasting of market operation rules and economic behaviors \cite{Graber2018}. In transportation, exploring mean-field game inverse problems can help optimize traffic flow control strategies, leading to improved efficiency and sustainability of road networks \cite{Cardaliaguet2018}. Moreover, mean-field game inverse problems find applications in financial risk assessment, energy management, environmental protection, and other fields \cite{Lachapelle2016,Yang2016,djehiche2017mean}. As research and technological advancements continue, this area will garner more attention and development. Furthermore, interdisciplinary collaborations will bring fresh perspectives and insights to the study of mean-field game inverse problems.

Specifically, there have been several numerical investigations\cite{Chow2022,Ding2022,Klibanov2024b,Mou2022} and theoretical analyses  \cite{Liu2023,Liu2022,Liu2023b,Liu2024, Liu2024a, Imanuvilov2023,Klibanov2023,Klibanov2024,Klibanov2023b, Yamamoto2009, Klibanov2024a,Klibanov1992,Klibanov2023d, Klibanov2023c,Immanuvilov1998,Liu2023c,Ren2023,Ren2024} focusing on MFG inverse problems using boundary/domain measurements. The current focus is primarily on using full boundary data or internal sub-region data to reconstruct unknown parameters or functions within the MFG model.  Studying the inversion of partial data presents both a realistic and challenging scenario. In many practical situations, obtaining complete boundary data can be difficult, costly, or even impossible. However, partial boundary data is often more readily available. Therefore, it is of utmost importance to develop effective methods and techniques that can accurately estimate the unknown parameters or functions within the mean-field games model using this limited information. The task of inferring the unknown parameters or functions in MFG models from partial boundary data is an active area of research. While there have been notable advancements in studying inverse problems using partial boundary data for MFG system, such as numerical studies \cite{Chow2022}, theoretical investigations \cite{Ding2023,Imanuvilov2023a} available for MFG, as well as inverse problems related to biomathematical coupled equations \cite{Ding2023a, Li2024, Liu2023a}, it is important to acknowledge that there is still much progress to be made in this field.

The present work is devoted to an inverse problem for mean-field game systems, with the aim of reconstructing unknown model components from observational data acquired only on a partial segment of the outer boundary. We first set out the geometric framework and basic notation used throughout the paper. Consider a bounded domain \(\Omega\subset\mathbb{R}^n\) with spatial dimension \(n\ge 2\), where \(\partial\Omega\) satisfies the \(C^{2+\alpha}\) regularity condition for some \(\alpha\in(0,1)\). Let \(D\Subset\Omega\) stand for a nonempty open subset with \(C^{2+\alpha}\) boundary such that \(\Omega\setminus\overline{D}\) forms a connected domain. We set the space‑time domain \(Q = (\Omega\setminus\overline{D})\times(0,T]\), and define its lateral boundaries \(\Sigma_1 = \partial D\times(0,T]\) and \(\Sigma_2 = \partial\Omega\times(0,T]\). The mean‑field game system we examine is a coupled forward‑backward parabolic system expressed as
\begin{equation}\label{um0}
\begin{cases}
-\dfrac{\partial v}{\partial t} - \nu\Delta v + \dfrac12\big|\nabla v\big|^2 = K(x,\rho),
& \text{in } Q,\\[4pt]
\dfrac{\partial \rho}{\partial t}-\nu\Delta \rho -\nabla\cdot\bigl(\rho\nabla v\bigr)=0,
& \text{in } Q,\\[4pt]
v=0,\quad \rho = f_0,
& \text{on } \Sigma_1,\\[4pt]
\mathcal{S} v = 0,\quad \mathcal{S} \rho = g_0,
& \text{on } \Sigma_2,\\[4pt]
\rho(\cdot,0)=\rho_0(x),\quad v(\cdot,T)=0,
& \text{in } \Omega\setminus \overline{D}.
\end{cases}
\end{equation}
Here $\Delta$ stands for the standard Laplacian operator and $\nabla\cdot$ the divergence operator, both acting on the spatial coordinate $x$. The first equation is a backward equation for the value function $v$, while the second is a forward equation governing the population density $\rho$. The boundary values $f_0$ and $g_0$ are given as nonnegative constants imposed on the lateral boundary. The constant $\nu$ is the diffusion coefficient governing the random dispersal of individual agents. The nonlinear term $K(x,\rho)$ serves as the running cost functional, which captures the interaction effect between each agent and the aggregate population distribution. The function $\rho_0(x)$ specifies the initial spatial profile of the agent population. The boundary observation operator $\mathcal{S}$ admits two canonical forms:

\[
\mathcal{S}v = v \quad \text{or} \quad \mathcal{S}v = \partial_\nu v,
\]
where $\partial_\nu v(x) = \nabla v(x) \cdot \nu(x)$ is the normal derivative of $u$, and $\nu$ denotes the unit outward normal vector to the boundary.

The mean‑field game coupled system \eqref{um0} couples the Hamilton‑Jacobi‑Bellman (HJB) and Kolmogorov‑Fokker‑Planck (KFP) equations to characterize system‑wide equilibrium behaviour. The HJB equation determines the value function for agents under stochastic optimal control, while the KFP equation governs the evolution of agent probability density on the game state space. For the MFG model studied here, we pay specific attention to the inner boundary \(\partial D\), which is also known as the Dirichlet interior boundary. This boundary condition characterizes the cost incurred by agents entering or exiting the game through $\partial D$ \cite{Ding2023, Ferreira2019}. For example, in a traffic model \cite{Cardaliaguet2018}, $D$ could represent a paid parking zone within a city $\Omega$. Vehicles entering or leaving this zone incur a fee. Conversely, the Neumann boundary condition on $\partial \Omega$ can represent a reflective boundary, preventing agents from leaving the domain \cite{Liu2023b}. By incorporating different types of boundaries, the MFG model aims to capture and simulate more realistic scenarios. This consideration of various boundary types is crucial for addressing inverse problems in complex situations, enabling predictions to be made in diverse settings.

For system \eqref{um0}, when the measurement mappings take the form $\mathcal{S}v=v$ and $\mathcal{S}\rho=\rho$, the corresponding measurement map is defined as:
\begin{equation}
\mathbf{{\Lambda}}^{\mathcal{D}}_{K,D}(\rho_0) =
\left( \int_{U_T}\frac{\partial v}{\partial \nu}(x,t)\, q(x,t) \, dxdt,\
\left. \frac{\partial \rho}{\partial \nu} \right|_{U_T} \right),
\end{equation}
with \(U_T = \Gamma\times(0,T]\) for \(\Gamma\subseteq\partial\Omega\), and $q$ representing a weight function. When $\mathcal{S}$ represents the normal derivative $$\mathcal{S} v = \partial_{\nu}v, \mathcal{S} \rho = \partial_{\nu}\rho,$$ the corresponding map is:
\[\mathbf{{\Lambda}}^{\mathcal{N}}_{K,D}(\rho_0) = \left(\int_{U_T}v(x,t)n(x,t)dxdt,\rho(x,t)\bigl|_{U_T}\right)\]
with \(n\) being another weight function. The measured data consist of two parts. The first involves the value function \(v\). In crowd dynamics \cite{Chow2022}, \(v\) can be interpreted as a travel cost; thus, \(\int_{U_T}\partial_\nu v q\) represents the weighted flux of this cost across a boundary segment over time, while \(\int_{U_T}v n\) represents a weighted average of the cost itself on the boundary. The second part involves the agent density \(\rho\), measured via its Neumann or Dirichlet trace on the partial boundary \cite{Ding2021}.

For the inverse problem aiming at simultaneous recovery of the running‑cost function $K(x,m)$ and boundary $\partial D$, its general formulation reads
\begin{align}\label{IP1}
\mathbf{{\Lambda}}^{\mathcal{D}}_{K,D}(\rho_0)\rightarrow K, D\ \   \mathrm{for}\  \rho_0\in \mathcal{H},
\end{align}
and
\begin{align}\label{IP2}
\mathbf{{\Lambda}}^{\mathcal{N}}_{K,D}(\rho_0)\rightarrow K, D\ \   \mathrm{for}\  \rho_0\in \mathcal{H}.
\end{align}
Here $\mathcal{H}$ represents the priori space, see Section \ref{sec2}.
This work focuses on the theoretical study of unique identifiability associated with the inverse‑problem formulations \eqref{IP1} and \eqref{IP2}. In order to settle this theoretical question, we construct sufficient assumptions within the general framework such that the problems \eqref{IP1} and \eqref{IP2} become uniquely identifiable,
\begin{align*}
\mathbf{{\Lambda}}^{\mathcal{D}}_{K_1,D_1}(\rho_0)=\mathbf{{\Lambda}}^{\mathcal{D}}_{K_2,D_2}(\rho_0)
\ \ \mathrm{if\ and\ only\ if}\  (K_1,D_1)=(K_2,D_2)\ \mathrm{for}\  \rho_0\in \mathcal{H},
\end{align*}
and
\begin{align*}
\mathbf{{\Lambda}}^{\mathcal{N}}_{K_1,D_1}(\rho_0)=\mathbf{{\Lambda}}^{\mathcal{N}}_{K_2,D_2}(\rho_0)
\ \ \mathrm{if\ and\ only\ if}\  (K_1,D_1)=(K_2,D_2)\ \mathrm{for}\  \rho_0\in \mathcal{H}.
\end{align*}
Here, for \(j=1,2\), \(D_j\) denotes a non‑vacuous open domain equipped with \(C^{2+\alpha}\) boundary regularity, and the complement region \(\Omega\setminus\overline{D_j}\) is assumed to be connected. Each \(K_j\) is taken from the a‑priori admissible family \(\mathcal{A}\) or \(\mathcal{B}\), introduced in Section \ref{sec2}.

For the sake of a general overview, we now sketch the principal outcomes of this article. Full technical arguments are presented throughout Sections \ref{sec3}--\ref{sec4}. In particular, we obtain uniqueness conclusions regarding the simultaneous recovery of the running cost and interior boundary for MFG systems subject to homogeneous as well as inhomogeneous boundary constraints.

\emph{(i) Simultaneous recovery results for homogeneous boundary MFG systems}
\begin{itemize}
\item{In the inverse problem of the homogeneous Dirichlet boundary system \eqref{um0}, it is possible to achieve a unique reconstruction of both the function $K(x,\rho)$, which belongs to an admissible class $\mathcal{A}$, and the internal topological boundary $\partial D$ by employing measurements of \eqref{IP1} with diverse inputs.}
\item{The inverse problem of the homogeneous mixed boundary system \eqref{um0} enables the unique reconstruction of the function $K(x,\rho)$, belonging to an admissible class $\mathcal{A}$, as well as the boundary $\partial D$, by employing measurements of \eqref{IP2} with diverse inputs.}
\end{itemize}

\emph{(ii) Simultaneous recovery results for inhomogeneous boundary MFG systems}
\begin{itemize}
\item{In the inverse problem related to the inhomogeneous Dirichlet system \eqref{um0}, the function $K(x,\rho)$, belonging to an admissible class $\mathcal{B}$, as well as the boundary $\partial D$, can be reconstructed uniquely by employing measurements of \eqref{IP1} with diverse inputs.}
\item{In the context of the inverse problem related to the inhomogeneous mixed boundary system \eqref{um0}, the function $K(x,\rho)$, which belongs to an admissible class $\mathcal{B}$, and the boundary $\partial D$ can be reconstructed uniquely by utilizing measurements of \eqref{IP2} with diverse inputs.}
\end{itemize}

\subsection{Technical developments and discussion}\label{sub1.2}

Within the framework of this inverse investigation, we exploit accessible mean‑flux measurements of \(v(x,t)\) over a partial boundary segment, together with Neumann or Dirichlet boundary traces associated with \(\rho(x,t)\) on the remaining boundary portion, to achieve joint identification of the running cost functional \(K(x,\rho)\) and the inner geometric boundary \(\partial D\) for the MFG system \eqref{um0}. Notably, the reconstruction of the boundary component \(\partial D\) represents one distinctive feature of the present work, which has remained absent from existing literature.

The study of the inverse problem is confronted with various challenges. These challenges encompass the nonlinear coupling between variables, the requirement for simultaneous reconstruction of shapes and parameters, the presence of nonnegative constraints, and the limitations of measurement information. For the inverse problem formulated from the nonlinear coupled system \eqref{um0}, the high‑order linearization method is adopted in our analysis. According to \cite{Lassas2020}, such a technique yields powerful capabilities when dealing with inverse problems for nonlinear PDEs. In this study, we aim to apply this method specifically to the nonlinear coupled system \eqref{um0}. This method is achieved by performing a linearization around a set of solutions $(v(x,t;\epsilon), \rho(x,t;\epsilon)) = (v(x,t;0), \rho(x,t;0))$ by injecting the following form of initial values
\begin{align}\label{m0}
\rho_0(x;\epsilon)=\rho_0(x;0)+\sum_{l=1}^N\epsilon_lg_l,\ \mathrm{for}\ x\in \Omega\setminus \overline{D},
\end{align}
where $(\epsilon_1,...\epsilon_N)\in \mathbb{R}^N$, $g_l\in C^{2+\alpha}(\Omega\setminus \overline{D})$ and $\rho_0(x;0)$ denotes the value of $\rho_0(x;\epsilon)$ when $\epsilon=0$. By considering the initial value input as \eqref{m0} and performing linearization of the solution $(v(x,t;0), \rho(x,t; 0))$, it can be concluded that  the nonnegativity of both $\rho_0$ and the boundary values $f_0$ and $g_0$ determines the nonnegativity of $\rho(x,t)$ within the domain $\overline{Q}$ for the system \eqref{um0}. For the present analytical framework, properly choosing the baseline solution pair \((v(x,t;0), \rho(x,t;0))\) is crucial when one implements the higher‑order linearization strategy corresponding to the nonlinear coupled system. This work examines two distinct expansion strategies to carry out system linearization. One strategy implements higher‑order linearization with respect to the baseline configuration \((v(x,t;0), \rho(x,t;0)) = (0, 0)\). By contrast, the other adopts linearization centred at \((v(x,t;0), \rho(x,t;0)) = (0, C)\), where $C$ represents some positive constant value. By employing these initial value selections, we can explore various scenarios and examine the influence of the starting point on the high-order linearization process of the nonlinear coupled system \eqref{um0}.

For the homogeneous‑boundary setting of system \eqref{um0}, i.e., when \(g_0=f_0=0\), we adopt the higher‑order linearization scheme expanded about the trivial baseline solution \((v(x,t;0), \rho(x,t;0))=(0,0)\). To guarantee non‑negativity for the population density \(\rho(x,t)\), the initial datum \(\rho_0\) is required to satisfy a non‑negativity constraint. We impose the condition that one element $g_l$ in \eqref{m0} of $\rho_0$ can be arbitrary, while the remaining elements must be strictly positive (see Section \ref{sec3}). In other words, this represents a form where the positive term dominates, while the impact on the arbitrary term is nearly negligible. Consequently, this enables us to linearize the solution at $(0,0)$, thereby simplifying the inherent complexity of the coupling problem. In references \cite{Liu2023b,Liu2023a}, the application of linearization near $(0,0)$ are discussed. These papers employ a higher-order variational method and employ distinct inputs for $\rho_0$ to address the inverse problem. However, it is important to emphasize that this method is not applicable for solving the inverse problem in MFG with a Dirichlet boundary condition. It is only suitable for scenarios involving Neumann boundaries and periodic boundaries, where full boundary measurements are available. For the strategy adopted herein, the arbitrariness of one entry in \(\rho_0\) enables us to employ Complex‑Geometric‑Optics (CGO) solutions for solving this inverse‑problem setting. Following the above methodology, one may achieve unique determination for the geometry of the interior boundary \(\partial D\) by carrying out first‑order linearization upon the equation satisfied by \(\rho\). To determine the behavior of the cost function $K$, a higher-order linearization is necessary.

For an inhomogeneous boundary in system \eqref{um0}, where $g_0$ and $f_0$ cannot both be zero simultaneously, we focus on two specific forms of boundary conditions in the inverse problem of the MFG system. First, we consider the Dirichlet boundary condition for $\mathcal{S}=\mathcal{I}$ in equation \eqref{um0}, where $\mathcal{I}$ represents a unit operator. In this case, we set $f_0=g_0$, ensuring that both $f_0$ and $g_0$ are positive on the boundary. The second form of boundary condition we examine is the Neumann boundary condition on the outer boundary $\partial \Omega$, where $f_0=0$, and a positive $g_0$ is imposed on the inner Dirichlet boundary $\partial D$. To address inverse‑problem questions for the MFG system under inhomogeneous boundary constraints, we make use of higher‑order linearization built around the reference state \((v, \rho)=(0,g_0)\). By approximating system dynamics in the neighbourhood of \((0,g_0)\), this linearization procedure alleviates difficulties originating from intricate nonlinear couplings and yields simplified coupled structures. So far, only \cite{Liu2022,Liu2024} has investigated the linearization at nonzero solutions. Since the linearization is performed on a positive solution $g_0$, adding a sufficiently small amount to this solution also guarantees a positive value for $\rho_0$ (see Section \ref{sec4}). This input provides a means to ensure the positivity of $\rho_0$ and facilitates the inversion process for added convenience. Within this approach, the shape of the internal boundary $\partial D$ can be uniquely determined through a first-order linearization of the equation for $\rho$. However, in order to comprehensively understand the behavior of the cost function $K$, a higher-order linearization becomes necessary.

In the remainder of this paper, we provide a structured presentation of our findings. The well‑posedness theory of the forward problem, together with its basic features, is presented in Section \ref{sec2}. Section \ref{sec3} then investigates inverse‑problem results on simultaneous reconstruction, placing particular emphasis on MFG systems subject to homogeneous boundary conditions. Advancing our study further, in Section \ref{sec4}, where we explore the unique results associated with inhomogeneous boundary MFG systems, shedding light on their distinctive characteristics.

\section{Well-posedness analysis for the MFG system}
\label{sec2}

Prior to investigating unique identifiability properties for the MFG system \eqref{um0}, we first derive local well‑posedness for its corresponding direct problem. Existence and uniqueness of solutions are established via an application of the implicit‑function theorem defined on complex Banach spaces; consult, for example, \cite{Gilbarg1977, Isakov1993} and the bibliographic sources cited therein. The analyticity property of the coupling term with respect to the population density variable occupies a central position within our analysis, and the exact set of hypotheses will be stated using the admissible function classes defined in the subsequent discussion.

We shall use the standard parabolic H\"{o}lder spaces. For \(\alpha\in(0,1)\), \(C^{k+\alpha}(\overline{\Omega})\) is the space of \(C^k(\overline{\Omega})\) functions with \(\alpha\)-Hölder continuous k-th derivatives. Its associated norm reads
\[
\|\phi\|_{C^{k+\alpha}(\overline{\Omega})}
=
\sum_{|l|\leq k}\|D^l\phi\|_{L^\infty(\Omega)}
+
\sum_{|l|=k}[D^l\phi]_{\alpha,\Omega},
\]
where
\[
[D^l\phi]_{\alpha,\Omega}
=
\sup_{x\neq y}
\frac{|D^l\phi(x)-D^l\phi(y)|}{|x-y|^\alpha},
\]
and
$D^l=\partial_{x_1}^{l_1}\cdots\partial_{x_n}^{l_n}$.

For space-time functions, we denote by
$C^{k+\alpha,\frac{k+\alpha}{2}}(\overline Q)$ the corresponding
parabolic H\"{o}lder space consisting of functions whose derivatives
$D^lD_t^j\phi$ with $|l|+2j\leq k$ are H\"{o}lder continuous with respect
to the parabolic scaling. The norm is given by
\[
\begin{aligned}
\bigl\|\phi\bigr\|_{C^{k+\alpha,\tfrac{k+\alpha}{2}}(\overline{Q})}
={}&
\sum_{\substack{|l|+2j = k}}
\sup_{\begin{smallmatrix}t\in(0,T]\\ x\neq y\end{smallmatrix}}
\frac{\bigl|D_x^l\phi(x,t)-D_x^l\phi(y,t)\bigr|}{|x-y|^\alpha}
\\[5pt]
&+
\sum_{\substack{|l|+2j = k}}
\sup_{\begin{smallmatrix}x\in\Omega\\ t\neq t'\end{smallmatrix}}
\frac{\bigl|D_x^l\phi(x,t)-D_x^l\phi(x,t')\bigr|}{|t-t'|^{\frac{\alpha}{2}}}
\\[5pt]
&+
\sum_{\substack{|l|+2j \le k}}
\bigl\| D_x^l D_t^j \phi \bigr\|_{L^\infty(Q)} < \infty .
\end{aligned}
\]
We next introduce the admissible classes of nonlinear coupling terms. All holomorphic expansions below are understood as expansions of
$C^\alpha(\Omega\setminus\overline D)$-valued holomorphic maps.

\begin{defn}[Admissible classes]\label{adm}
Let $K:\Omega\setminus\overline{D}\times\mathbb{C}\to\mathbb{C}$.
\begin{itemize}
\item[(i)] We say that $K\in\mathcal{A}$ whenever the mapping $z\mapsto K(\cdot,z)$ defines a holomorphic function taking values in $C^\alpha(\Omega\setminus\overline{D})$, with $K(x,0)=0$ holding for every $x\in\Omega\setminus\overline{D}$.

\item[(ii)] Let $g_0>0$ be a given constant. We write $K\in\mathcal{B}(g_0)$ provided the mapping $z\mapsto K(\cdot,z)$ is holomorphic and maps into $C^\alpha(\Omega\setminus\overline{D})$. Furthermore, the conditions $K(x,g_0)=0$ and $\partial_z K(x,g_0)=0$ are satisfied for all $x$.
\end{itemize}
\end{defn}

For any element $K\in\mathcal{A}$, the above analyticity property guarantees the Taylor series representation
\[
K(x,z)=\sum_{i=1}^{\infty} K^{(i)}(x)\,\frac{z^i}{i!},
\qquad
K^{(i)}(x) = \left.\frac{\partial^i K}{\partial z^i}\right|_{z=0} \;\in C^\alpha(\Omega\setminus\overline{D}).
\]
Similarly, for $K\in\mathcal B(g_0)$, the first-order vanishing
conditions imply that the expansion around $g_0$ starts from the
second-order term:
\[
K(x,z)=\sum_{i=2}^{\infty} K^{(i)}(x)\,\frac{(z-g_0)^{i}}{i!},
\qquad
K^{(i)}(x)=\left.\frac{\partial^{i}K}{\partial z^{i}}\right|_{z=g_0}(x,g_0)
\in C^{\alpha}\bigl(\Omega\setminus\overline{D}\bigr).
\]
The two admissible classes correspond to different reference
configurations in the inverse problem and will be used separately in
the identifiability analysis. Next, we present a result on the local well-posedness of the MFG system \eqref{um0} with both inhomogeneous and homogeneous boundaries.

\begin{thm}\label{thm1}
Assume $K\in\mathcal{A}$ and $f_0=g_0=0$ in \eqref{um0}. Set
\[
X = C^{2+\alpha}(\overline{\Omega}\setminus D), \qquad
Y = C^{2+\alpha,1+\alpha/2}(\overline{Q}).
\]
For \((m,w)\in Y\times Y\), define \(\|(m,w)\|_{Y\times Y}:=\|m\|_Y+\|w\|_Y\).
Then there are constants \(\delta>0\), \(C>0\) such that the following properties hold:
\begin{enumerate}
\item[(i)] For every $\rho_0$ in the ball
\[
B_\delta := \{ \phi\in X : \|\phi\|_X \le \delta \},
\]
system \eqref{um0} admits a unique solution $(v,\rho)\in Y\times Y$ satisfying
\[
\|(v,\rho)\|_{Y\times Y} \le C \|\rho_0\|_X .
\]
Consequently, the solution belongs to
\[
\mathcal{X}_\delta := \{ (m,w)\in Y\times Y : \|(m,w)\|_{Y\times Y} \le C\delta \}.
\]
\item[(ii)] The solution mapping
\[
S_1 : B_\delta \longrightarrow \mathcal{X}_\delta, \qquad S_1(\rho_0) := (v,\rho),
\]
is holomorphic.
\end{enumerate}
\end{thm}

\begin{thm}\label{thm2}
Assume $F\in\mathcal{B}$ and that in \eqref{um0} either
\begin{itemize}
    \item[(a)] $\mathcal{S}=\mathcal{I}$ with $f_0=g_0$, or
    \item[(b)] $\mathcal{S}=\partial_\nu$ with $f_0=0$ and $g_0>0$.
\end{itemize}
Let $X = C^{2+\alpha}\bigl(\overline{\Omega}\setminus D\bigr)$, $Y = C^{2+\alpha,1+\tfrac{\alpha}{2}}\bigl(\overline{Q}\bigr)$ and $W = C^{2+\alpha,1+\tfrac{\alpha}{2}}\bigl(\overline{\Sigma_1\cup\Sigma_2}\bigr)$,
and for $(m,w)\in Y\times Y$ define $\|(m,w)\|_{Y\times Y}:=\|m\|_Y+\|w\|_Y$.
Then there exist constants $\delta>0$ and $C>0$ such that the following hold:
\begin{enumerate}
\item[(i)]Given \((\rho_0,g_0)\in\)
\(B_\delta^2:=\bigl\{ (\phi,\psi)\in X\times W : \phi=\psi \ \text{on}\ \partial\Omega\cup\partial D,\; \|\phi\|_X+\|\psi\|_W\le\delta \bigr\},\)
the system \eqref{um0} has a unique solution \((v,\rho)\in Y\times Y\) satisfying
\[
\|(v,\rho)\|_{Y\times Y}\le C\bigl(\|\rho_0\|_X+\|g_0\|_W\bigr).
\]
Consequently, the solution belongs to
\[
\mathcal{X}_\delta:=\bigl\{ (v,w)\in Y\times Y : \|(v,w)\|_{Y\times Y}\le C\delta \bigr\}.
\]
\item[(ii)] The solution operator
\[
S_2: B_\delta^2\longrightarrow\mathcal{X}_\delta,\qquad S_2(\rho_0,g_0):=(v,\rho),
\]
is holomorphic.
\end{enumerate}
\end{thm}

\begin{proof}
We prove the case \(\mathcal{S}=\mathcal{I}\), \(f_0=g_0\);
the case \(\mathcal{S}=\partial_\nu\), \(f_0=0\), \(g_0>0\) follows from similar arguments. For brevity, we define the Banach spaces
\[
X = C^{2+\alpha}\bigl(\overline{\Omega}\setminus D\bigr),
\qquad
Y = C^{2+\alpha,1+\tfrac{\alpha}{2}}\bigl(\overline{Q}\bigr),
\qquad
Z = C^{2+\alpha,1+\tfrac{\alpha}{2}}\bigl(\overline{\Sigma_1\cup\Sigma_2}\bigr),
\qquad
E = C^{\alpha,\tfrac{\alpha}{2}}\bigl(\overline{Q}\bigr).
\]
Using these abbreviations, define
\[
\begin{aligned}
Y_1 &:= \{\, \phi\in X : \phi = g_0 \ \text{on}\ \partial\Omega\cup\partial D \,\},\\[2pt]
Y_2 &:= \{\, \psi\in Z : \psi = g_0 \ \text{on}\ \overline{\Sigma_1}\cup\overline{\Sigma_2} \,\},\\[2pt]
Y_3 &:= \{\, (v,w)\in Y\times Y : v(\cdot,T)=0 \ \text{in}\ \overline{\Omega}\setminus D,\;
            v=0 \ \text{on}\ \Sigma_1\cup\Sigma_2 \,\},\\[2pt]
Y_4 &:= Y_1\times Y_2\times E\times E .
\end{aligned}
\]
A nonlinear mapping $\mathfrak{F}\colon Y_1\times Y_2\times Y_3\to Y_4$ is defined componentwise as follows:
\begin{equation}
\label{eq:operator-F-mfg}
\mathfrak{F}\big(\tilde{\rho}_0,\tilde{g}_0,\tilde{v},\tilde{\rho}\big)
:=
\begin{pmatrix}
\tilde{\rho}(\cdot,0) - \tilde{\rho}_0 \\[4pt]
\bigl(\tilde{\rho}-\tilde{g}_0\bigr)\Big|_{\Sigma_1\cup\Sigma_2} \\[4pt]
-\partial_t \tilde{v} - \nu\,\Delta\,\tilde{v}
+\frac12 \bigl|\nabla \tilde{v}\bigr|^{2}
-K(x,\tilde{\rho}) \\[4pt]
\partial_t \tilde{\rho} - \nu\,\Delta\,\tilde{\rho}
-\nabla\cdot\Bigl(\tilde{\rho}\,\nabla \tilde{v}\Bigl)
\end{pmatrix}.
\end{equation}
Pointwise multiplication preserves H\"{o}lder regularity, so both $|\nabla \tilde{v}|^2$ and $\nabla\cdot(\tilde{\rho}\nabla \tilde{v})$ lie in the space $E$. The only remaining step is to confirm $K(x,\tilde{\rho})\in E$.
Since $K$ is holomorphic with respect to its second argument $z$, we have the estimate
\[
\|K^{(i)}\|_{C^\alpha(\overline{\Omega}\setminus\overline{D})}
\leq \frac{i!}{R^i}\sup_{|z|=R}\|K(\cdot,z)\|_{C^\alpha(\overline{\Omega}\setminus\overline{D})}
\quad\text{for each } R>0.
\]
From the above inequality, one can find a positive constant $L>0$ for which
\[
\biggl\|\frac{K^{(i)}}{i!}\,\tilde{\rho}^{\,i}\biggr\|_E
\le \frac{L^i}{R^i}\,\|\tilde{\rho}\|_E^{\,i}\,
\sup_{|z|=R}\|K(\cdot,z)\|_{C^{\alpha}(\Omega\setminus\overline{D})}.
\]
After substituting \(R = 2L\|\tilde{\rho}\|_E\), we obtain convergence in $E$ for the power series
\(\sum_{i=2}^{\infty}\frac{K^{(i)}(x)}{i!}\tilde{\rho}^{\,i}\).
As a consequence, \(K(x,\tilde{\rho})\) belongs to $E$, a fact that secures the well‑posedness of the mapping \(\mathfrak{F}\). Given the local boundedness property of \(\mathfrak{F}\), our task reduces to establishing weak holomorphic. Fix two points
\[
(\rho_0,\tilde{g}_0,\tilde{v},\tilde{\rho}),\,
(\bar{\rho}_0,\bar{g}_0,\bar{v},\bar{\rho})\in Y_1\times Y_2\times Y_3
\]
and consider
\[
\beta \longmapsto
\mathfrak{F}\bigl((\rho_0,\tilde{g}_0,\tilde{v},\tilde{\rho})
+ \beta(\bar{\rho}_0,\bar{g}_0,\bar{v},\bar{\rho})\bigr).
\]
The only non-trivial part is the dependence on $\tilde{\rho}$ through $K$.
Because the expansion
\[
\sum_{i=2}^{\infty}\frac{K^{(i)}(x)}{i!}\,(\rho_0+\beta\bar{\rho}_0)^i
\]
converges in $E$ locally uniformly in $\beta\in\mathbb{C}$,
the map $\beta\mapsto K(x,\rho_0+\beta\bar{\rho}_0)$ is holomorphic with values in $E$.
Thus $\mathfrak{F}$ is holomorphic.

A direct computation shows $\mathfrak{F}(g_0,g_0,0,g_0)=0$.
The Fr\'{e}chet derivative with respect to $(v,\rho)$ at this point acts as
\[
\begin{aligned}
\partial_{(v,\rho)}\mathfrak{F}(g_0,g_0,0,g_0)(v,\rho)
= \bigl(
&\rho(x,0),\;
\rho|_{\Sigma_1\cup\Sigma_2},\\
&-\partial_t v-\nu\Delta v,\;
\partial_t \rho-\nu\Delta \rho - g_0\Delta v
\bigr).
\end{aligned}
\]
If this derivative vanishes, then $v=0$ and subsequently $\rho=0$;
hence it is injective. For surjectivity, let $(r,s),\,(r',s')\in X\times E$.
We can find $a\in Y$ solving the backward parabolic problem
\[
-\partial_t a-\nu\Delta a = s,\quad a|_{\Sigma_1\cup\Sigma_2}=0,\quad a(\cdot,T)=r,
\]
and then $b\in Y$ solving
\[
\partial_t b-\nu\Delta b - g_0\Delta a = s',\quad b|_{\Sigma_1\cup\Sigma_2}=0,\quad b(\cdot,0)=r'.
\]
The existence of such constructed pairs $a$ and $b$ establishes surjectivity for the linearized operator; accordingly, this operator constitutes a linear isomorphism. By virtue of the implicit function theorem, there exists a uniquely determined holomorphic map
\(S_2 : B_\delta \to Y_3\) fulfilling
\[
\mathfrak{F}\bigl(\rho_0,g_0,S_2(\rho_0,g_0)\bigr)=0
\qquad\text{for all } (\rho_0,g_0)\in B_\delta.
\]
Setting $(v,\rho)=S_2(\rho_0,g_0)$ yields the unique solution of \eqref{um0}.
Finally, write $(v_0,\rho_0)=S_2(g_0,g_0)$ and use the Lipschitz continuity of $S_2$ to estimate
\[
\begin{aligned}
\|(v,\rho)\|_{Y\times Y}
&\le C\bigl( \|\rho_0\|_X + \|g_0\|_Z + \|v_0\|_Y + \|\rho_0\|_Y \bigr)\\[4pt]
&\le C\bigl( \|\rho_0\|_X + \|g_0\|_Z \bigr),
\end{aligned}
\]
where the last step follows from the bound already established for $S_2$.
This completes the proof.
\end{proof}

\section{High order linearization}\label{sec:linearzation_discussion}

\subsection{Linearization in the space of probability measures}

This section starts with a recapitulation of the variational setting defined on probability‑measure spaces; this formalism serves for carrying out linearization procedures associated with mean‑field game systems. Although the linearized MFG equations have been extensively investigated on
the torus (see, e.g., \cite{Cardaliaguet2019}), analogous arguments remain
valid in bounded domains; see \cite{Ricciardi2022}. In the present work,
the measure variable is represented by a sufficiently smooth density, but
the linearization is naturally formulated through derivatives on the space of probability measures.

Suppose \(\tilde{\Omega}\subset\mathbb R^n\) is a bounded domain. Let \(\mathcal P(\tilde{\Omega})\) be the collection of Borel probability measures on \(\tilde{\Omega}\). Given the functional
\(U\colon \mathcal P(\tilde{\Omega})\rightarrow\mathbb R,\) we recall the notion of its first variation.

\begin{defn}\label{def_der_1}
A functional $U$ is called differentiable over $\mathcal P(\tilde{\Omega})$ provided that one can find a continuous mapping
\[
\mathcal{K}\colon \mathcal P(\tilde{\Omega})\times\tilde{\Omega}\to\mathbb R
\]
satisfying
\begin{equation}\label{derivation}
\lim_{s\to 0^+}
\frac{U\bigl(m_1+s(m_2-m_1)\bigr)-U(m_1)}{s}
=
\int_{\tilde{\Omega}} \mathcal{K}(m_1,x)\,d(m_2-m_1)(x)
\end{equation}
for all $m_1,m_2\in\mathcal P(\tilde{\Omega})$.
\end{defn}

Since the right-hand side of \eqref{derivation} is invariant under the
addition of constants to $K$, the derivative is defined up to a constant.
Throughout this paper we choose the normalized representative satisfying
\[
\int_{\tilde{\Omega}}\mathcal{K}(m,x)\,dm(x)=0.
\]
This representative is denoted by
\[
\frac{\delta U}{\delta m}(m,x).
\]
Consequently, for $m_1,m_2\in\mathcal P(\tilde{\Omega})$,
\[
U(m_2)-U(m_1)
=
\int_0^1
\int_{\tilde{\Omega}}
\frac{\delta U}{\delta m}
(m_1+s(m_2-m_1),x)
\,d(m_2-m_1)(x)\,ds .
\]
Higher-order derivatives with respect to the measure variable are defined
recursively; we refer to \cite{Ricciardi2022} for further details. The natural metric for perturbations of probability measures is the
$1$-Wasserstein distance.

\begin{defn}\label{W_distance}
Given two elements $\mu_1$ and $\mu_2$ from $\mathcal P(\tilde{\Omega})$, their Wasserstein distance is introduced as
\[
\mathbb W_1(\mu_1,\mu_2)
:=
\sup\left\{
\int_{\tilde{\Omega}}\varphi(x)\,d(\mu_1-\mu_2)(x)
\,\bigg|\, \mathrm{Lip}(\varphi)\le 1
\right\},
\]
in which
\[
\mathrm{Lip}(\varphi)
\triangleq
\sup_{\substack{x,y\in\tilde{\Omega}\\ x\neq y}}
\frac{|\varphi(x)-\varphi(y)|}{|x-y|}
\]
stands for the minimal Lipschitz constant associated with $\varphi$.
\end{defn}

\begin{rem}\label{remark_measure_density}
In the above definitions, the variable $m$ is regarded as a probability
measure. In the MFG system studied in this paper, however, the same notation
is used for the corresponding density function. Throughout the analysis, this
density is assumed to have the H\"{o}lder regularity required in Sections
\ref{sec2} and \ref{sec3}.

The connection between the measure-theoretic linearization and the classical
linearized PDE system can be described as follows. Let
$(v_1,\rho_1)$ and $(v_2,\rho_2)$ be two solutions associated with initial densities $\rho_0^1$ and $\rho_0^2$, respectively. The first variation around
$(v_2,\rho_2)$ leads to a linearized pair $(\theta,\eta)$ whose initial
perturbation is given by $\rho_0^1-\rho_0^2$. Under suitable regularity assumptions,
the nonlinear increment
\[
(v_1-v_2,\rho_1-\rho_2)
\]
can be approximated by $(\theta,\eta)$, with the remainder controlled by the
distance between the initial measures. In our setting, since the densities
belong to H\"{o}lder spaces rather than merely to $\mathcal P(\tilde{\Omega})$, the required stability estimates follow from
the well-posedness theory established in Theorems~\ref{thm1} and
\ref{thm2}.
\end{rem}

\subsection{Higher-order linearization and an auxiliary spectral lemma}

We investigate the linearized structure of the forward system by expanding the
solution around the zero state. Throughout this subsection, we take
$\mathcal{S}=\mathcal{I}$ and $f_0=g_0=0$. The same argument can be adapted to other admissible background states. Let the initial density be perturbed by the finite-dimensional family
\begin{equation}\label{m00}
\rho_0(x;\epsilon)
=
\sum_{l=1}^{N}\epsilon_l g_l ,
\qquad
g_l\in C^{2+\alpha}(\overline{\Omega}\setminus D).
\end{equation}
Let the perturbation parameter vector be $\epsilon = (\epsilon_1,\dots,\epsilon_N)\in\mathbb{R}^N$, with its $\ell^1$-norm defined by
\[
|\epsilon| := \sum_{l=1}^N |\epsilon_l|,
\]
and suppose $|\epsilon|$ is sufficiently small. The associated solution pair $(v,\rho)$ to the perturbed MFG system takes the form
\begin{equation}\label{ume}
\begin{cases}
\begin{aligned}
-\dfrac{\partial v}{\partial t} - \nu\,\Delta v + \dfrac12 \bigl|\nabla v\bigr|^{2}
&= K(x,\rho),
&& \text{in } Q,\\[5pt]
\dfrac{\partial \rho}{\partial t} - \nu\,\Delta \rho - \nabla\cdot\Bigl(\rho\,\nabla v\Bigr)
&= 0,
&& \text{in } Q,\\[5pt]
v = 0,\quad \rho &= 0,
&& \text{on } \Sigma_1\cup\Sigma_2,\\[5pt]
v\Big|_{t=T}=0,\quad \rho\Big|_{t=0}&=\sum_{l=1}^{N}\epsilon_l\,g_l
&& \text{in } \Omega\setminus\overline{D}.
\end{aligned}
\end{cases}
\end{equation}
By virtue of the local well-posedness established in Section 2, system \eqref{ume} possesses a unique solution for sufficiently small $|\epsilon|$, belonging to the regularity class $(v,\rho).$
Furthermore, the admissibility requirement $K(x,0)=0$ ensures that vanishing perturbation $\epsilon=0$ recovers the trivial zero solution of the system. From Thereom\ref{thm1}(ii) imply that $S_1$ is holomorphic near the
origin. Hence its Fr\'{e}chet derivative
$A_1=S_1'(0)$
exists and satisfies
\[
S_1(\rho_0)
=
A_1(\rho_0)
+
o\left(
\|\rho_0\|_{C^{2+\alpha}(\overline{\Omega}\setminus D)}
\right)~~~as~~ \rho_0\rightarrow0.
\]
To obtain the first variation, we choose the perturbation direction
$\rho_0=\epsilon_1g_1$. The linearity of $A_1$ gives
\[
A_1(\epsilon_1g_1)
=
\epsilon_1(v^{(1)},\rho^{(1)}),
\]
where
\[
(v^{(1)},\rho^{(1)})
=
A_1(g_1)
=
\left.
(\partial_{\epsilon_1}v,
\partial_{\epsilon_1}\rho)
\right|_{\epsilon=0}.
\]
Taking the first-order variation of system \eqref{ume} with respect to $\epsilon_1$, evaluated at the zero perturbation $\epsilon=0$, gives rise to the following linearized system:
\begin{equation}\label{eq3.8}
\begin{cases}
\begin{aligned}
-\dfrac{\partial v^{(1)}}{\partial t} - \nu\,\Delta v^{(1)}
&= K^{(1)}(x)\, \rho^{(1)},
&& \text{in } Q,\\[5pt]
\dfrac{\partial \rho^{(1)}}{\partial t} - \nu\,\Delta \rho^{(1)}
&= 0,
&& \text{in } Q,\\[5pt]
v^{(1)} = 0,\quad \rho^{(1)} &= 0,
&& \text{on } \Sigma_1\cup\Sigma_2,\\[5pt]
v^{(1)}\Big|_{t=T}=0,\quad \rho^{(1)}\Big|_{t=0}&=g_1
&& \text{in } \Omega\setminus\overline{D}.
\end{aligned}
\end{cases}
\end{equation}
More generally, for every index $l=1,\dots,N$, the first-order directional derivative pair
\[
(v^{(l)}, \rho^{(l)}) = \left.(\partial_{\epsilon_l} v,\, \partial_{\epsilon_l} \rho)\right|_{\epsilon=0}
\]
obeys an identical linearized system, with initial datum given by $g_l$.

Higher-order linearised equations are constructed by taking successive mixed partial derivatives with respect to the perturbation parameters. As a representative example, the second-order mixed variation reads
\[
(v^{(1,2)}, \rho^{(1,2)}) = \left.(\partial_{\epsilon_2}\partial_{\epsilon_1} v,\, \partial_{\epsilon_2}\partial_{\epsilon_1} \rho)\right|_{\epsilon=0}.
\]
solves
\begin{equation}\label{umN2N}
\left\{
\begin{aligned}
-\dfrac{\partial v^{(1,2)}}{\partial t} - \nu\,\Delta v^{(1,2)}
+\Bigl(\nabla v^{(1)}\Bigr) \cdot \Bigl(\nabla v^{(2)}\Bigr)
&= K^{(1)}\, \rho^{(1,2)} + K^{(2)}\, \rho^{(1)} \rho^{(2)},
&& \forall (x,t)\in Q,\\[5pt]
\dfrac{\partial \rho^{(1,2)}}{\partial t} - \nu\,\Delta \rho^{(1,2)}
-\nabla\cdot\Bigl(\rho^{(1)}\nabla v^{(2)}\Bigr)
&-\nabla\cdot\Bigl(\rho^{(2)}\nabla v^{(1)}\Bigr)
= 0,
&& \forall (x,t)\in Q,\\[5pt]
v^{(1,2)} = 0,\quad \rho^{(1,2)} = 0,
&
&& \text{on } \Sigma_1\cup\Sigma_2,\\[5pt]
v^{(1,2)}\Big|_{t=T} = 0,\quad \rho^{(1,2)}\Big|_{t=0} = 0
&
&& x\in \Omega\setminus\overline{D}.
\end{aligned}
\right.
\end{equation}
The boundary and initial conditions are obtained by differentiating those of
\eqref{ume}. Higher-order mixed derivatives are defined in the same manner and generate a sequence of linear parabolic systems.

We next introduce a spectral lemma that provides the special solutions used in the uniqueness argument.

\begin{lem}\label{lem1}
This parabolic system admits a sequence of solutions to
\begin{equation}\label{eigD}
\left\{
\begin{aligned}
\dfrac{\partial Z(x,t)}{\partial t} - \nu\,\Delta Z(x,t)
&= 0,
&& \forall (x,t)\in Q,\\[5pt]
Z(x,t) &= 0,
&& \text{on } \Sigma_1,\\[5pt]
\mathcal{S}Z(x,t) &= 0,
&& \text{on } \Sigma_2.
\end{aligned}
\right.
\end{equation}
admitting separable solutions of the form
\[
Z_k(x,t) = e^{-\mu_k\, t}\, y_k(x),
\]
where $\mu_k\in\mathbb R$ and
$y_k\in C^{2+\alpha}(\overline{\Omega}\setminus D)$.
\end{lem}

\begin{proof}
Consider the elliptic operator associated with the mixed boundary value
problem
\[
-\Delta y=\lambda y
\]
in $\Omega\setminus\overline D$ together with the boundary conditions
\[
y=0\quad\text{on }\partial D,
\qquad
\mathcal{S}y=0\quad\text{on }\partial\Omega .
\]
The corresponding elliptic operator possesses a discrete spectrum
$\{\lambda_k\}_{k\ge1}$ with eigenfunctions $y_k$ by standard spectral
theory.  Taking
\[
\mu_k=\nu\lambda_k
\]
and setting
\[
Z_k(x,t)=e^{-\mu_k t}y_k(x),
\]
a direct substitution verifies that $Z_k$ satisfies the parabolic problem
\eqref{eigD}.
\end{proof}

\section{Simultaneous recovery results related to homogeneous boundary MFG systems}\label{sec3}

The purpose of this section is to establish the identifiability of the
unknown coupling function and the inaccessible boundary in the homogeneous
MFG system. More precisely, we show that the observation operators associated
with the exterior boundary measurements uniquely determine the pair
$(K,D)$ when the nonlinear interaction term belongs to the admissible class
$\mathcal A$. The arguments will be developed separately for the Dirichlet
and Neumann type observations.

\subsection{Uniqueness for the external Dirichlet boundary MFG system}

We begin our analysis with measurement data collected from exterior Dirichlet boundary observations. The theorem below gives a rigorous characterization of uniqueness for the associated inverse problem.

\begin{thm}\label{Dthm1}
For $j=1,2$, let $(K_j,D_j)$ be a pair of admissible configurations. Each cost function $K_j$ belongs to the admissible class $\mathcal{A}$, and each domain $D_j\Subset\Omega$ is a nonempty open set endowed with $C^{2+\alpha}$ boundary regularity. The complement $\Omega\setminus\overline{D_j}$ is connected. Let
$q\in C_0^{2+\alpha}(S_T)$
be nonnegative and not identically zero. Denote by
$\mathbf{\Lambda}_{K_j,D_j}^{\mathcal D}$
the Dirichlet measurement operator associated with system
\eqref{um0} under the homogeneous condition
$\mathcal S=\mathcal I$ and
$f_0=g_0=0$. Assume that
\begin{equation}\label{codD1}
\mathbf{\Lambda}_{K_1,D_1}^{\mathcal D}(\rho_0)
=
\mathbf{\Lambda}_{K_2,D_2}^{\mathcal D}(\rho_0)
\end{equation}
for all admissible initial data $m_0$ in the common domain of definition.
Then
\[
D_1=D_2,
\qquad
K_1(x,\rho)=K_2(x,\rho).
\]
\end{thm}

\begin{proof}
\emph{Step I. Unique determination of the interior domain $D$.} Our argument begins with the construction of a parameterized family of initial conditions:
\begin{equation}\label{fg}
\rho_0(x)=\sum_{l=1}^{N}\epsilon_l g_l,\qquad x\in\Omega\setminus\overline{D},
\end{equation}
with \(g_l\in C^{2+\alpha}(\Omega\setminus\overline{D})\) for \(l=1,\dots,N\), and \(g_l>0\) for all \(l\neq 2\). The parameter vector $\epsilon=(\epsilon_1,\dots,\epsilon_N)\in\mathbb{R}_+^N$ is chosen sufficiently small so that
\[
\Bigl\|\sum_{l=1}^N\epsilon_l g_l\Bigr\|_{C^{2+\alpha}(\Omega\setminus\overline{D})}<\delta.
\]
No sign restriction is imposed on $g_2$, which can be selected arbitrarily; we only require $\epsilon_2\ll\epsilon_l$ for every $l\neq2$ to ensure the strict positivity $m_0>0$ on the whole domain $\Omega\setminus\overline{D}$. From the local well‑posedness conclusion proven within Section~\ref{sec2}, it follows that for any admissible \(\epsilon\), system \eqref{um0} has exactly one classical solution pair satisfying
\(\bigl(v(\cdot,\cdot;\epsilon),\,\rho(\cdot,\cdot;\epsilon)\bigr)\in C^{2+\alpha,1+\alpha/2}(\overline{Q})\times C^{2+\alpha,1+\alpha/2}(\overline{Q}).\)
In particular, when $\epsilon=0$ the system reduces to the trivial solution $$(v(\cdot,\cdot;0),\rho(\cdot,\cdot;0))=(0,0)$$ over $Q$. Performing a first-order linearization of system \eqref{ume} with respect to $\epsilon_1$ and evaluating at $\epsilon=0$, we derive for $j=1,2$ the linearized systems:
\begin{equation}\label{mD11D}
\begin{cases}
\begin{aligned}
-\dfrac{\partial v_j^{(1)}}{\partial t} - \nu\,\Delta v_j^{(1)}
&= K^{(1)} \rho_j^{(1)},
&& \text{in } \bigl(\Omega\setminus\overline{D_j}\bigr)\times(0,T],\\[5pt]
\dfrac{\partial \rho_j^{(1)}}{\partial t} - \nu\,\Delta \rho_j^{(1)}
&= 0,
&& \text{in } \bigl(\Omega\setminus\overline{D_j}\bigr)\times(0,T],\\[5pt]
v_j^{(1)} = 0,\quad \rho_j^{(1)} &= 0,
&& \text{on } \bigl(\partial D_j\cup\partial\Omega\bigr)\times(0,T],\\[5pt]
v_j^{(1)}\Big|_{t=T}=0,\quad &\rho_j^{(1)}\Big|_{t=0}=g_1>0
&& \text{in } \Omega\setminus\overline{D_j}.
\end{aligned}
\end{cases}
\end{equation}
 Let $G$ designate the connected component of $\Omega\setminus(\overline{D_1}\cup\overline{D_2})$ whose boundary contains the outer boundary $\partial\Omega$, and define the difference function $V:=\rho_1^{(1)}-\rho_2^{(1)}$. Combining the equality of the measurement maps given in \eqref{codD1} with the fact that the Neumann traces of $\rho$ coincide on $\Gamma$, we infer that $V$ satisfies
 \begin{equation}\label{mD21}
\begin{cases}
\begin{aligned}
\partial_t V - \nu\,\Delta\, V
&= 0,
&& \text{in } G\times(0,T],\\[4pt]
\partial_\nu V &= 0,
&& \text{on } S_T,\\[4pt]
V &= 0,
&& \text{on } \Sigma_2,\\[4pt]
V(x,0) &= 0
&& \text{in } G.
\end{aligned}
\end{cases}
\end{equation}
Appealing to the unique continuation principle for linear parabolic equations (cf. \cite{Saut1987}), we deduce that $V$ vanishes identically on $\overline{G}\times(0,T]$, which implies $\rho_1^{(1)}=\rho_2^{(1)}$ on $\overline{G}\times(0,T]$. In conjunction with the boundary conditions in \eqref{mD11D}, this yields $\rho_1^{(1)}=\rho_2^{(1)}=0$ along $\partial(D_2\setminus D_1)\times(0,T]$. Since the initial datum $g_1$ is strictly positive, the strong parabolic maximum principle entails $\rho_1^{(1)}>0$ everywhere in $(\Omega\setminus\overline{D_1})\times(0,T]$. This stands in contradiction with the vanishing of $\rho_1^{(1)}$ on $\partial(D_2\setminus D_1)\times(0,T]$. Therefore we must have $D_1=D_2$.

\medskip

\emph{Step II. Identification of the first-order Taylor coefficient $F^{(1)}$.} With the identity $D_1=D_2=:D$ now in hand, we proceed to carry out a first-order linearization with respect to the parameter $\epsilon_2$. Define
\[
v^{(2)}:=\partial_{\epsilon_2}v|_{\epsilon=0},\qquad \rho^{(2)}:=\partial_{\epsilon_2}\rho|_{\epsilon=0},
\]
and observe immediately that $\rho^{(2)}$ does not depend on the index $j$, as the evolution equation for $\rho$ contains no terms involving $K$. For each $j=1,2$, the linearized system reads
\begin{equation}\label{umD11}
\begin{cases}
\begin{aligned}
-\dfrac{\partial v_j^{(2)}}{\partial t} - \nu\,\Delta v_j^{(2)}
&= K_j^{(1)} \rho^{(2)},
&& \text{within } Q,\\[5pt]
\dfrac{\partial \rho^{(2)}}{\partial t} - \nu\,\Delta \rho^{(2)}
&= 0,
&& \text{within } Q,\\[5pt]
v_j^{(2)} = 0,\quad \rho^{(2)} &= 0,
&& \text{over } \Sigma_1\cup\Sigma_2,\\[5pt]
v_j^{(2)}\Big|_{t=T}=0,\quad \rho^{(2)}\Big|_{t=0}&=g_2
&& \text{in } \Omega\setminus\overline{D}.
\end{aligned}
\end{cases}
\end{equation}
where $g_2$ remains arbitrary. Subtracting the first equation in \eqref{umD11} for $j=1$ from that for $j=2$, we arrive at a boundary value problem for the difference $w:=v_1^{(2)}-v_2^{(2)}$:
\begin{equation}\label{wk11}
\begin{cases}
\begin{aligned}
-\dfrac{\partial w}{\partial t} - \nu\,\Delta w
&= \bigl(K_1^{(1)}-K_2^{(1)}\bigr)\, \rho^{(2)},
&& \text{within } Q,\\[5pt]
w &= 0,
&& \text{over } \Sigma_1\cup\Sigma_2,\\[5pt]
w\Big|_{t=T}&=0,
&& \text{in } \Omega\setminus\overline{D}.
\end{aligned}
\end{cases}
\end{equation}
We now introduce an auxiliary adjoint problem for a test function $b$:
\begin{equation}\label{z}
\left\{
\begin{aligned}
\partial_{t} b - \nu \Delta b &= 0, && \text{in } Q,\\[5pt]
b &= 0, && \text{on } \Sigma_1,\\[5pt]
b &= q_0, && \text{on } \Sigma_2,\\[5pt]
b(x,0) &= 0, && \text{in } \Omega\setminus\overline{D}.
\end{aligned}
\right.
\end{equation}
where $q_0$ agrees with $q$ on $U_T$ and vanishes identically on $(\partial\Omega\setminus\Gamma)\times(0,T]$. Multiplying both sides of \eqref{wk11} by $b$, integrating over the space-time domain $Q$, and applying Green's formula together with the measurement equality
\[
\int_{U_T}\partial_\nu v_1^{(2)}q \,dSdt = \int_{U_T}\partial_\nu v_2^{(2)}q \,dSdt,
\]
led to the integral relation
\[
\iint_{(\Omega\setminus\overline{D}\times[0,T])}
\big(K_1^{(1)}(x)-K_2^{(1)}(x)\big)\,\rho^{(2)}(x,t)\,b(x,t)\,dxdt
= 0.
\]
In light of Lemma~\ref{lem1}, set $\rho^{(2)}(x,t)=e^{-\mu t}\beta(x,\mu)$, with $\beta(\cdot,\mu)$ denoting an arbitrary Dirichlet eigenfunction of $-\Delta$ on $\Omega\setminus\overline{D}$. Substituting this ansatz into the foregoing identity gives
\[
\bigl(K_1^{(1)}(x)-K_2^{(1)}(x)\bigr) \int_0^T e^{-\mu t} b(x,t)\,dt = 0.
\]
The strong maximum principle guarantees $b>0$ throughout $Q$, so the time integral is strictly positive. It follows immediately that $K_1^{(1)}(x)=K_2^{(1)}(x)$ for all $x\in\Omega\setminus\overline{D}$.

\medskip

\emph{Step III. Recovery of the second-order Taylor coefficient $F^{(2)}$.}
Once the equalities $D_1 = D_2$ and $K_1^{(1)} = K_2^{(1)}$ have been verified, we move on to analyze the mixed second-order partial derivative $\partial_{\epsilon_1}\partial_{\epsilon_2}$ evaluated at $\epsilon = 0$.
A direct observation shows that $\rho^{(1,2)}$ is independent of the index $j$.
By carrying out the second-order linearization, we arrive at the following coupled system, which holds for each $j = 1,2$,
\begin{equation}\label{umD2}
\left\{
\begin{aligned}
-\dfrac{\partial v_j^{(1,2)}}{\partial t}-\nu\,\Delta v_j^{(1,2)}
+\Bigl(\nabla v^{(1)}\Bigr)\cdot \Bigl(\nabla v^{(2)}\Bigr)
&= K^{(1)} \rho^{(1,2)} + K_j^{(2)}\, \rho^{(1)} \rho^{(2)},
&& \text{in } Q,\\[5pt]
\dfrac{\partial \rho^{(1,2)}}{\partial t}-\nu\,\Delta \rho^{(1,2)}
-\nabla\cdot\Bigl(\rho^{(1)}\nabla v^{(2)}\Bigr)
&-\nabla\cdot\Bigl(\rho^{(2)}\nabla v^{(1)}\Bigr)
= 0,
&& \text{in } Q,\\[5pt]
v_j^{(1,2)} &= 0,\qquad \rho^{(1,2)} = 0,
&& \text{on } \Sigma_1\cup\Sigma_2,\\[5pt]
v_j^{(1,2)}\Big|_{t=T} &= 0,\quad \rho^{(1,2)}\Big|_{t=0}=0,
&& \text{in }\Omega\setminus\overline{D}.
\end{aligned}
\right.
\end{equation}
The quantities \(v^{(1)}\), \(v^{(2)}\), \(\rho^{(1)}\), \(\rho^{(2)}\) and \(\rho^{(1,2)}\) get their unique solutions from the lower‑order linearized problems; their values bear no dependence upon the coefficient \(K_j^{(2)}\). We now subtract the first equation for $j=2$ from the corresponding equation for $j=1$, and introduce the difference function $z := v_1^{(1,2)} - v_2^{(1,2)}$.
This new unknown satisfies the boundary value problem below:
\begin{equation}\label{FD211}
\left\{
\begin{aligned}
-\dfrac{\partial z}{\partial t} - \nu\,\Delta z
&= \big(K_1^{(2)} - K_2^{(2)}\big)\, \rho^{(1)} \rho^{(2)},
&& \forall\,(x,t)\in Q,\\[5pt]
z &= 0,
&& (x,t)\in\Sigma_1 \cup \Sigma_2,\\[5pt]
z\Big|_{t=T} &= 0,
&& \forall\,x\in \Omega \setminus \overline{D}.
\end{aligned}
\right.
\end{equation}
In analogy with Step II, multiplication of equation \eqref{FD211} by the adjoint function $b$ introduced in \eqref{z}, followed by integration over $Q$ and a further application of Green's formula combined with the measurement identity \eqref{codD1}, gives
\[
\iint_{(0,T)\times(\Omega\setminus\overline{D})}
\big(K_1^{(2)}(x)-K_2^{(2)}(x)\big)\,
\rho^{(1)}(x,t)\,\rho^{(2)}(x,t)\,b(x,t)\,dx\,dt = 0.
\]
With $\rho^{(2)}(x,t)=e^{-\mu t}\beta(x,\mu)$ chosen as before, and with the strict positivity $\rho^{(1)}>0$ and $b>0$ in $Q$ guaranteed by the parabolic maximum principle, it follows that $K_1^{(2)}(x)=K_2^{(2)}(x)$ for every $x\in\Omega\setminus\overline{D}$.

\medskip

\emph{Step IV. Inductive recovery of higher-order coefficients $K^{(N)}$ for $N\geq 3$.} We proceed by mathematical induction. Let $N\ge3$ be an integer, and assume as induction hypothesis that the identities $K_1^{(i)}=K_2^{(i)}$ hold for all orders $i=1,\dots,N-1$. Taking the $N$-th order mixed derivative $\partial_{\epsilon_1}\cdots\partial_{\epsilon_N}$ at $\epsilon=0$ and subtracting the resulting equations for $j=1$ and $j=2$, we arrive at an integral identity of the form
\[
\iint_{(0,T)\times(\Omega\setminus\overline{D})}
\big(K_1^{(N)}(x)-K_2^{(N)}(x)\big)\,
\rho^{(1)}\rho^{(2)}\rho^{(3)}\cdots \rho^{(N)}\,b(x,t)\,dxdt=0
\]
where $\rho^{(2)}$ is once again chosen as $e^{-\mu t}\beta(x,\mu)$. Each of the remaining factors $\rho^{(1)},\rho^{(3)},\dots,\rho^{(N)}$ is strictly positive over the whole space-time domain $Q$. Since $b>0$ in $Q$ by the maximum principle, the integrand has fixed sign unless $K_1^{(N)}(x)=K_2^{(N)}(x)$ pointwise on $\Omega\setminus\overline{D}$. By induction, all Taylor coefficients of $K_1$ and $K_2$ coincide, which implies the global identity $K_1(x,\rho)=K_2(x,\rho)$ on $(\Omega\setminus \overline{D})\times\mathbb{R}^+$. The proof is complete.
\end{proof}

\subsection{Uniqueness for the external Neumann boundary MFG system}

A second uniqueness result is established below, which asserts that the measurement map $\mathbf{\Lambda}_{K,D}^{\mathcal{N}}$ introduced in \eqref{IP2} uniquely identifies both the cost function $K(x,\rho)$ and the interior obstacle boundary $\partial D$.

\begin{thm}\label{Nthm2}
For $j=1,2$, let $K_j\in\mathcal{A}$, and let $D_j\Subset\Omega$ stand for non‑empty open sets endowed with $C^{2+\alpha}$ boundary regularity, whose exterior domains $\Omega\setminus\overline{D_j}$ are connected. The weight function $g\in C_0^{2+\alpha}(S_T)$ is taken to be nonnegative and does not vanish identically. Denote by $\Lambda_{K_j,D_j}^{\mathcal{N}}$ the measurement map for the MFG system \eqref{um0} with $\mathcal{S}=\partial_\nu$ and homogeneous boundary data $f_0=g_0=0$.
If
\begin{equation}\label{cod2N}
\Lambda_{K_1,D_1}^{\mathcal{N}}(\rho_0)=\Lambda_{K_2,D_2}^{\mathcal{N}}(\rho_0)
\end{equation}
holds for every initial datum $\rho_0\in C^{2+\alpha}(\Omega\setminus\overline{D})$,
then
\[
D_1=D_2,
\qquad
K_1(x,\rho)=K_2(x,\rho).
\]
\end{thm}

\begin{proof}
\emph{Step I}. We first determine the interior boundary. Consider the perturbation \eqref{fg} and linearize the MFG system at the trivial solution
\((v,\rho)=(0,0)\). For $j=1,2$, the resulting systems for $(v_j^{(1)},\rho_j^{(1)})$ read
\begin{equation}\label{mN1N}
\left\{
\begin{aligned}
-\partial_{t} v_j^{(1)}-\nu \Delta v_j^{(1)}
&= K^{(1)}\, \rho_j^{(1)},
&& (x,t)\in \bigl(\Omega\setminus\overline{D_j}\bigr)\times(0,T],\\[5pt]
\partial_{t} \rho_j^{(1)}-\nu \Delta \rho_j^{(1)}
&= 0,
&& (x,t)\in \bigl(\Omega\setminus\overline{D_j}\bigr)\times(0,T],\\[5pt]
v_j^{(1)} = \rho_j^{(1)} &= 0,
&& (x,t)\in \partial D_j\times(0,T],\\[5pt]
\partial_{\nu} v_j^{(1)} = \partial_{\nu} \rho_j^{(1)}& = 0,
&& (x,t)\in \Sigma_2,\\[5pt]
v_j^{(1)}\big|_{t=T} = 0,\quad \rho_j^{(1)}\big|_{t=0}&=g_1>0,
&& x\in \Omega\setminus\overline{D_j}.
\end{aligned}
\right.
\end{equation}
Assume, for contradiction, that $D_1\neq D_2.$ Let $G$ stand for that connected component of \(\Omega\setminus(\overline{D_1}\cup\overline{D_2})\) adjacent to the outer boundary. Define
\begin{equation}\label{eq_V_aux}
\left\{
\begin{aligned}
\dfrac{\partial V}{\partial t} - \nu\,\Delta V &= 0,
&& (x,t)\in G\times(0,T],\\[5pt]
V &= 0,
&& (x,t)\in S_T,\\[5pt]
\partial_{\nu} V &= 0,
&& (x,t)\in \Sigma_2,\\[5pt]
V\Big|_{t=0} &= 0,
&& x\in G.
\end{aligned}
\right.
\end{equation}
The unique continuation property for parabolic equations gives
$V\equiv0\quad\text{in }\overline G\times(0,T].$
Hence,
\[
\rho_1^{(1)}=\rho_2^{(1)}
\quad\text{in }\overline G\times(0,T].
\]
Consequently, \(\rho_1^{(1)}\) vanishes on
\(\partial(D_2\setminus D_1)\times(0,T]\).
On the other hand, the positivity of the initial datum \(g_1\) and
strong maximum principle imply
\[
\rho_1^{(1)}>0
\quad\text{in }
(\Omega\setminus\overline{D_1})\times(0,T],
\]
which is a contradiction. Therefore, $D_1=D_2 .$

\emph{Step II}. Set $D_1=D_2=:D.$ The first-order linearization with respect to the second perturbation parameter gives $(v_j^{(2)},\rho^{(2)})$ for $j=1,2,$ where \(\rho^{(2)}\) is independent of \(j\). The difference $w=v_1^{(2)}-v_2^{(2)}$ satisfies
\begin{equation}\label{wk1N}
\left\{
\begin{aligned}
-\dfrac{\partial w}{\partial t} - \nu\,\Delta w
&= \big(K_1^{(1)}-K_2^{(1)}\big)\, \rho^{(2)},
&& (x,t)\in Q,\\[5pt]
w &= 0,
&& (x,t)\in \Sigma_1,\\[5pt]
\partial_{\nu} w &= 0,
&& (x,t)\in\Sigma_2,\\[5pt]
w\Big|_{t=T} &= 0,
&& x\in \Omega \setminus \overline{D}.
\end{aligned}
\right.
\end{equation}
Let $\tilde{b}$ be the solution of the corresponding adjoint homogeneous
problem with boundary datum chosen according to the observation weight.
Multiplying \eqref{wk1N} by $\tilde{b}$, integrating over the space-time
cylinder, and using Green's identity together with the equality of the
measurement operators, we obtain
\[
\iint_{(0,T)\times(\Omega\setminus\overline{D})}
\big(K_2^{(1)}(x)-K_1^{(1)}(x)\big)\,
\rho^{(2)}(x,t)\,\tilde{b}(x,t)\,dxdt=0
\]
Choosing the separated solutions
\[
\rho^{(2)}(x,t)=e^{-\mu t}\alpha(x,\rho)
\]
provided by Lemma~\ref{lem1}, where
\(\alpha(\cdot,\rho)\) is a mixed eigenfunction of the Laplacian, gives
\[
(K_2^{(1)}(x)-K_1^{(1)}(x))
\int_0^T e^{-\mu t}\tilde{b}(x,t)\,dt=0 .
\]
Since \(\tilde{b}>0\) by maximum principle, we conclude that $K_1^{(1)}=K_2^{(1)}
\text{in }\Omega\setminus\overline D .$

\emph{Step III.} Now fix $D_1=D_2$ and $K_1^{(1)}=K_2^{(1)}$. Taking the mixed second-order derivative with respect to the perturbation parameters, we obtain the system for $(v_j^{(1,2)},\rho^{(1,2)}),$
\begin{equation}\label{umN2N2}
\left\{
\begin{aligned}
-\dfrac{\partial v_j^{(1,2)}}{\partial t}-\nu\,\Delta v_j^{(1,2)}
+\Bigl(\nabla v^{(1)}\Bigr)\cdot \Bigl(\nabla v^{(2)}\Bigr)
&= K^{(1)}\, \rho^{(1,2)} + K_j^{(2)}\, \rho^{(1)} \rho^{(2)},
&& \forall\,(x,t)\in Q,\\[5pt]
\dfrac{\partial \rho^{(1,2)}}{\partial t}-\nu\,\Delta \rho^{(1,2)}
-\nabla\cdot\Bigl(\rho^{(1)}\nabla v^{(2)}\Bigr)
&-\nabla\cdot\Bigl(\rho^{(2)}\nabla v^{(1)}\Bigr)
= 0,
&& \forall\,(x,t)\in Q,\\[5pt]
v_j^{(1,2)} &= 0,\qquad \rho^{(1,2)} = 0,
&& \forall\,(x,t)\in\Sigma_1,\\[5pt]
\partial_{\nu} v_j^{(1,2)} &= 0,\quad \partial_{\nu} \rho^{(1,2)} = 0,
&& \forall\,(x,t)\in \Sigma_2,\\[5pt]
v_j^{(1,2)}\Big|_{t=T} &= 0,\quad \rho^{(1,2)}\Big|_{t=0}=0,
&& \forall\,x\in \Omega\setminus\overline{D}.
\end{aligned}
\right.
\end{equation}
Note that none of \(v^{(1)}\), \(v^{(2)}\), \(\rho^{(1)}\), \(\rho^{(2)}\), \(\rho^{(1,2)}\) depends on the coefficient \(K_j^{(2)}\). Subtracting the equations corresponding to \(j=1,2\), the difference $z=v_1^{(1,2)}-v_2^{(1,2)}$
satisfies
\begin{equation}\label{F2N}
\left\{
\begin{aligned}
-\dfrac{\partial z}{\partial t} - \nu\,\Delta z
&= \big(K_1^{(2)}-K_2^{(2)}\big)\, \rho^{(1)} \rho^{(2)},
&& \forall\,(x,t)\in Q,\\[5pt]
z &= 0,
&& \forall\,(x,t)\in \Sigma_1,\\[5pt]
\partial_{\nu} z &= 0,
&& \forall\,(x,t)\in\Sigma_2,\\[5pt]
z\Big|_{t=T} &= 0,
&& \forall\,x\in \Omega \setminus \overline{D}.
\end{aligned}
\right.
\end{equation}
Repeating the duality argument used in Step II gives
\[
\iint_{(0,T)\times(\Omega\setminus\overline{D})}
\big(K_1^{(2)}-K_2^{(2)}\big)\,
\rho^{(1)}\,\rho^{(2)}\,\tilde{b}\,dxdt=0
\]
The positivity of \(\rho^{(1)}\), \(\rho^{(2)}\), and \(\tilde{b}\) yields $K_1^{(2)}=K_2^{(2)}\text{in }\Omega\setminus\overline D .
$

\emph{Step IV.} In order to establish the validity of the result in equation \eqref{fg} for $N\geq 3$, it is necessary to compute the higher-order Taylor coefficients as follows:
\begin{align*}
v^{(1,2,...,N)}:&=\partial_{\epsilon_1}\partial_{\epsilon_2}...\partial_{\epsilon_N}v|_{\epsilon=0},\\
\rho^{(1,2,...,N)}:&=\partial_{\epsilon_1}\partial_{\epsilon_2}...\partial_{\epsilon_N}\rho|_{\epsilon=0}.
\end{align*}
By utilizing the powerful technique of mathematical induction, we can successfully derive the desired result. In other words, \(K_1^{(i)}(x)=K_2^{(i)}(x)\) holds in \(\Omega\setminus\overline{D}\) for every \(i\in N\). As a consequence, one may conclude that \(K_1(x,\rho)=K_2(x,\rho)\) on \((\Omega\setminus\overline{D})\times\mathbb{R}^+\). This completes the proof.

\end{proof}

\section{Simultaneous recovery results related to inhomogeneous boundary MFG systems}\label{sec4}

\subsection{Uniqueness for the external Dirichlet boundary MFG systems}

This subsection is devoted to the first uniqueness theorem for the mean field game system subject to inhomogeneous boundary conditions. Specifically, the Dirichlet type measurement map $\mathbf{\Lambda}_{K,D}^{\mathcal{D}}$ is shown to enable the simultaneous unique recovery of both the running cost function $K(x,\rho)$ and the geometric profile of the interior obstacle boundary $\partial D$.

\begin{thm}\label{thm:inhom-Dirichlet}
For $j=1,2$, let $K_j\in\mathcal{A}$, and let $D_j\Subset\Omega$ denote non‑empty open domains equipped with $C^{2+\alpha}$ boundary regularity, for which the exterior complement $\Omega\setminus\overline{D_j}$ is connected. The weight function $q\in C_0^{2+\alpha}(S_T)$ is taken to be nonnegative and not identically vanishing. Let $\mathbf{\Lambda}_{K_j,D_j}^{\mathcal{D}}$ stand for the boundary measurement map corresponding to the MFG system \eqref{um0} with observation operator $\mathcal{S}=\mathcal{I}$, considered under inhomogeneous boundary data satisfying $f_0=g_0>0$. If
\begin{equation}\label{cod1}
\Lambda_{K_1,D_1}^{\mathcal{D}}(\rho_0)=\Lambda_{K_2,D_2}^{\mathcal{D}}(\rho_0)
\end{equation}
holds for all admissible initial data $m_0\in C^{2+\alpha}(\Omega\setminus\overline{D_j})$, it follows that
\[
K_1 = K_2 \quad\text{and}\quad D_1 = D_2.
\]
\end{thm}

\begin{proof}
\emph{Step I. Determination of the interior boundary $\partial D$.}
We introduce initial data of the form
\begin{equation}\label{fg2}
\rho_0(x)=g_0+\sum_{l=1}^{N}\epsilon_l g_l,\qquad x\in\Omega\setminus\overline{D},
\end{equation}
in which \(g_l\) belong to \(C^{2+\alpha}(\Omega\setminus\overline{D})\) with indices \(l=1,\dots,N\), and the vector \(\epsilon=(\epsilon_1,\dots,\epsilon_N)\in\mathbb{R}^N\) is chosen small enough so as to ensure that
\[
\|g_0\|_{C^{2+\alpha,1+\alpha/2}(\Sigma_1\cup\Sigma_2)}+\Big\|\sum_{l=1}^N\epsilon_l g_l\Big\|_{C^{2+\alpha}(\Omega\setminus D)}
\]
is small enough to guarantee $\rho_0>0$ on $\Omega\setminus\overline{D}$. With $f_0=g_0>0$ and $\mathcal{S}=\mathcal{I}$, system \eqref{um0} becomes
\begin{equation}\label{um0h}
\left\{
\begin{aligned}
-\dfrac{\partial v}{\partial t} - \nu\,\Delta v + \tfrac12 \lvert\nabla v\rvert^2
&= K(x,\rho),
&& \forall\,(x,t)\in Q,\\[5pt]
\dfrac{\partial \rho}{\partial t} - \nu\,\Delta \rho - \nabla\cdot\Bigl(\rho\nabla v\Bigr)
&= 0,
&& \forall\,(x,t)\in Q,\\[5pt]
v = 0,\quad \rho &= g_0,
&& \text{on } \Sigma_1\cup\Sigma_2,\\[5pt]
v\Big|_{t=T} = 0,\quad \rho\Big|_{t=0}&=g_0+\sum_{l=1}^{N}\epsilon_l\, g_l,
&& \forall\,x\in \Omega\setminus\overline{D}.
\end{aligned}
\right.
\end{equation}
From well‑posedness arguments, the unique solution pair
\((v(\cdot,\cdot;\epsilon),\rho(\cdot,\cdot;\epsilon)) \in C^{2+\alpha,1+\alpha/2}(\overline{Q})\times C^{2+\alpha,1+\alpha/2}(\overline{Q})\)
is guaranteed to exist for all sufficiently small \(\epsilon\). The baseline solution at \(\epsilon=0\) satisfies the relation \((v(\cdot,\cdot;0),\rho(\cdot,\cdot;0))=(0,g_0)\) within the domain $Q$.
Let \(g_1\) be a positive function defined on \(\Omega\setminus\overline{D}\). The first-order linear derivative \(\partial_{\epsilon_1}\) evaluated at \(\epsilon=0\) about the baseline state \((0,g_0)\) generates the corresponding systems for \(j=1,2\),
\begin{equation}\label{mD11}
\left\{
\begin{aligned}
-\dfrac{\partial v_j^{(1)}}{\partial t}-\nu\,\Delta v_j^{(1)}
&= 0,
&& \forall\,(x,t)\in \bigl(\Omega\setminus\overline{D_j}\bigr)\times(0,T],\\[5pt]
\dfrac{\partial \rho_j^{(1)}}{\partial t}-\nu\,\Delta \rho_j^{(1)}-g_0\,\Delta v_j^{(1)}
&= 0,
&& \forall\,(x,t)\in \bigl(\Omega\setminus\overline{D_j}\bigr)\times(0,T],\\[5pt]
v_j^{(1)} = 0,\quad \rho_j^{(1)} &= 0,
&& \text{on } \bigl(\partial D_j\cup\partial\Omega\bigr)\times(0,T],\\[5pt]
v_j^{(1)}\Big|_{t=T} = 0,\quad \rho_j^{(1)}\Big|_{t=0}&=g_1>0,
&& \forall\,x\in \Omega\setminus\overline{D_j}.
\end{aligned}
\right.
\end{equation}
From the first equation and the homogeneous boundary conditions we obtain $v_j^{(1)}\equiv0$ in $\overline{Q}$, hence $\Delta v_j^{(1)}=0$. Consequently \eqref{mD11} simplifies to
\begin{equation}\label{mD1}
\left\{
\begin{aligned}
\dfrac{\partial \rho_j^{(1)}}{\partial t}-\nu\,\Delta \rho_j^{(1)}
&= 0,
&& \forall\,(x,t)\in \bigl(\Omega\setminus\overline{D_j}\bigr)\times(0,T],\\[5pt]
\rho_j^{(1)} &= 0,
&& \text{on } \bigl(\partial D_j\cup\partial\Omega\bigr)\times(0,T],\\[5pt]
\rho_j^{(1)}\Big|_{t=0} &= g_1>0,
&& \forall\,x\in \Omega\setminus\overline{D_j}.
\end{aligned}
\right.
\end{equation}
Denote by $G$ the connected component of $\Omega\setminus(\overline{D_1}\cup\overline{D_2})$ having $\partial\Omega$ as part of its boundary, and set $V:=\rho_1^{(1)}-\rho_2^{(1)}$. From the equality of the measurement maps \eqref{cod1} we deduce that $V$ satisfies problem \eqref{mD21}. The unique continuation principle for linear parabolic equations yields $V\equiv0$ in $\overline{G}\times(0,T]$; thus $\rho_1^{(1)}=\rho_2^{(1)}$ on $\overline{G}\times(0,T]$. By \eqref{mD1}, we obtain $\rho^{(1)}_1=\rho^{(1)}_2=0$ on $\partial(D_2\setminus D_1)\times (0,T]$. As $g_1>0$, the maximum principle forces $\rho^{(1)}_1>0$ in $\overline{Q}$, in contradiction to $\rho_1^{(1)}=0$ on $\partial(D_2\setminus D_1)\times (0,T]$. Therefore, we have $D_1=D_2$.

\emph{Step II. Identification of the second-order Taylor coefficient $F^{(2)}$.}
With the identity $D_1 = D_2 =: D$ already verified in the preceding step, we proceed to perform the second-order linearization of the system.
We introduce the second-order variational derivatives
$$v_j^{(1,2)} \triangleq \partial_{\epsilon_1}\partial_{\epsilon_2} v \big|_{\epsilon=0},\qquad
\rho^{(1,2)} \triangleq \partial_{\epsilon_1}\partial_{\epsilon_2} \rho \big|_{\epsilon=0},$$
where the quantity $\rho^{(1,2)}$ is independent of the index $j$.
Expanding the governing equations up to second order in $\epsilon$ leads to the following coupled system for each $j = 1,2$,
\begin{equation}\label{umDc2}
\left\{
\begin{aligned}
-\dfrac{\partial v_j^{(1,2)}}{\partial t}
-\nu\,\Delta v_j^{(1,2)}
+\bigl(\nabla v^{(1)}\bigr)\cdot \bigl(\nabla v^{(2)}\bigr)
&= K_j^{(2)} \, \rho^{(1)}\, \rho^{(2)},
&& \forall\,(x,t)\in Q,\\[5pt]
\dfrac{\partial \rho^{(1,2)}}{\partial t}
-\nu\,\Delta \rho^{(1,2)}
-\nabla\cdot\Bigl(\rho^{(1)}\nabla v^{(2)}\Bigr)
&-\nabla\cdot\Bigl(\rho^{(2)}\nabla v^{(1)}\Bigr)
-g_0\,\Delta v_j^{(1,2)}
= 0,
&& \forall\,(x,t)\in Q,\\[5pt]
v_j^{(1,2)} = 0,\quad \rho^{(1,2)}& = 0,
&& \text{on } \Sigma_1 \cup \Sigma_2,\\[5pt]
v_j^{(1,2)}\Big|_{t=T}=0,\quad
\rho^{(1,2)}\Big|_{t=0}&=0,
&& \forall\,x\in \Omega\setminus\overline{D}.
\end{aligned}
\right.
\end{equation}
Recall from Step I that $v^{(1)} \equiv 0$ and $v^{(2)} \equiv 0$ hold on $\overline{Q}$, so that $\nabla v^{(1)} = \nabla v^{(2)} = 0$ identically.
Inserting these relations into \eqref{umDc2} reduces the system to
\begin{equation}\label{umDcsimple}
\left\{
\begin{aligned}
-\dfrac{\partial v_j^{(1,2)}}{\partial t}
-\nu\,\Delta v_j^{(1,2)}
&= K_j^{(2)}\, \rho^{(1)}\, \rho^{(2)},
&& \forall\,(x,t)\in Q,\\[5pt]
\dfrac{\partial \rho^{(1,2)}}{\partial t}
-\nu\,\Delta \rho^{(1,2)}
-g_0\,\Delta v_j^{(1,2)}
&= 0,
&& \forall\,(x,t)\in Q,\\[5pt]
v_j^{(1,2)} = 0,\quad \rho^{(1,2)}& = 0,
&& \text{on } \Sigma_1 \cup \Sigma_2,\\[5pt]
v_j^{(1,2)}\Big|_{t=T}=0,\quad
\rho^{(1,2)}\Big|_{t=0} &= 0,
&& \forall\,x\in \Omega\setminus\overline{D}.
\end{aligned}
\right.
\end{equation}
Subtracting the first equations for $j=1,2$ gives
\begin{equation}\label{FD2}
\left\{
\begin{aligned}
-\dfrac{\partial z}{\partial t} - \nu\,\Delta z
&= \big(K_1^{(2)}-K_2^{(2)}\big)\, \rho^{(1)} \rho^{(2)},
&& (x,t)\in Q,\\[5pt]
z &= 0,
&& \text{on } \Sigma_1\cup\Sigma_2,\\[5pt]
z\Big|_{t=T} &= 0,
&& x\in \Omega\setminus\overline{D},
\end{aligned}
\right.
\end{equation}
where we have set $z := v_1^{(1,2)} - v_2^{(1,2)}$.
Denote $b$ as the solution satisfying the adjoint boundary‑value problem \eqref{z}. Upon multiplying equation \eqref{FD2} by b and performing integration across the space‑time domain $Q$, and invoking Green's identity combined with the measurement condition
$\int_{U_T} \partial_\nu v_1^{(1,2)} q \,dS = \int_{U_T} \partial_\nu v_2^{(1,2)} q \,dS$,
one arrives at
\begin{equation*}
\iint_{(0,T)\times(\Omega\setminus\overline{D})}
\big(K_1^{(2)}(x)-K_2^{(2)}(x)\big)\,
\rho^{(1)}(x,t)\,\rho^{(2)}(x,t)\,b(x,t)\,dxdt=0
\end{equation*}
By Lemma~\ref{lem1}, we may choose $\rho^{(2)}(x,t) = e^{-\mu t} \beta(x,\mu)$,
where $\beta(\cdot,\mu)$ stands for an arbitrary Dirichlet eigenfunction of $-\Delta$ on $\Omega\setminus\overline{D}$.
Since the family of Dirichlet eigenfunctions is complete in $L^2(\Omega\setminus\overline{D})$, it follows that
\begin{equation*}
\bigl(K_1^{(2)}(x) - K_2^{(2)}(x)\bigr)
\int_0^T e^{-\mu t} \rho^{(1)}(x,t) b(x,t) \,dt = 0.
\end{equation*}
An application of the strong maximum principle yields strict positivity for b together with \(\rho^{(1)}\) at every point inside the space‑time cylinder $Q$. This immediately yields $K_1^{(2)}(x) = K_2^{(2)}(x)$ for every $x \in \Omega\setminus\overline{D}$.

\emph{Step III. Identification of higher-order Taylor coefficients $K^{(N)}$ for $N \geq 3$.}
To extend identity \eqref{fg} to the case $N \geq 3$, we carry out the same construction for higher-order variational derivatives.
Proceeding inductively, we derive the integral relation
\begin{equation*}
\iint_{(0,T)\times(\Omega\setminus\overline{D})}
\big(K_1^{(N)}(x)-K_2^{(N)}(x)\big)\,
\rho^{(1)}\,\rho^{(2)}\cdots \rho^{(N)}\,b(x,t)\,dxdt=0
\end{equation*}
Repeating the reasoning from Step II, we take $\rho^{(2)} = e^{-\mu t} \beta(x,\mu)$,
while $\rho^{(1)}, \rho^{(3)}, \dots, \rho^{(N)}$ are all positive solutions on $\overline{Q}$.
The same completeness argument then gives $K_1^{(N)}(x) = K_2^{(N)}(x)$ on $\Omega\setminus\overline{D}$.
By analyticity of the Taylor expansion, we conclude $K_1(x,\rho) = K_2(x,\rho)$ on $(\Omega\setminus\overline{D}) \times \mathbb{R}^+$, which finishes the proof.

\end{proof}

\subsection{Uniqueness for the external Neumann boundary MFG system}
A second uniqueness result demonstrates that both the cost functional \(K(x,\rho)\) and the inner boundary \(\partial D\) admit unique reconstruction from the measurement map \(\mathbf{\Lambda}_{K,D}^{\mathcal{N}}\) defined in \eqref{IP2}.

\begin{thm}\label{thm:neumann-inhom}
For $j=1,2$, let $K_j\in\mathcal{B}$, and let $D_j\Subset\Omega$ denote non-empty open subsets equipped with $C^{2+\alpha}$ boundary regularity, whose complementary region $\Omega\setminus\overline{D_j}$ is connected. Let $n\in C_0^{2+\alpha}(S_T)$ be a nonnegative weight function that is not identically zero. Denote by $\mathbf{\Lambda}_{K_j,D_j}^{\mathcal{N}}$ the measurement map associated with the MFG system \eqref{um0} for $\mathcal{S}=\partial_\nu$ and the boundary data $f_0=0$, $g_0>0$. If for every initial datum $\rho_0\in C^{2+\alpha}(\Omega)$ there holds
\begin{equation}\label{cod2}
\Lambda_{K_1,D_1}^{\mathcal{N}}(\rho_0)=\Lambda_{K_2,D_2}^{\mathcal{N}}(\rho_0),
\end{equation}
then
\[
K_1=K_2\quad\text{and}\quad D_1=D_2.
\]
\end{thm}

\begin{proof}
The argument proceeds in several steps, using higher-order linearisations around the non-zero olution $(v,\rho)=(0,g_0)$. $\emph{Step\ I}.$ We choose initial data of the form
\begin{equation}\label{fg2b}
\rho_0(x)=g_0+\sum_{l=1}^{N}\epsilon_l g_l(x),\qquad x\in\Omega\setminus\overline{D},
\end{equation}
with $g_l\in C^{2+\alpha}(\Omega\setminus\overline{D})$ and $\epsilon=(\epsilon_1,\dots,\epsilon_N)\in\mathbb{R}^N$ so small that $\rho_0>0$ on $\Omega\setminus\overline{D}$. For $\epsilon=0$ the unique solution of \eqref{um0} (with $\mathcal{S}=\partial_\nu$, $f_0=0$, $g_0>0$) is $(v(\cdot,\cdot;0),\rho(\cdot,\cdot;0))=(0,g_0)$ in $Q$. Let \(g_1>0\) be given on \(\Omega\setminus\overline{D}\), and consider the first‑order derivatives
\(v_j^{(1)}:=\partial_{\epsilon_j}v\Big|_{\epsilon=0},\  \rho_j^{(1)}:=\partial_{\epsilon_j}\rho\Big|_{\epsilon=0}\) for $j=1,2$.
Performing a first-order linearization of system (1.1) about the reference state $(0,g_0)$, we derive the following linearized system for each domain $D_j$:
\begin{equation}\label{mN121}
\left\{
\begin{aligned}
-\dfrac{\partial v_j^{(1)}}{\partial t}-\nu\,\Delta v_j^{(1)}
&= 0,
&& (x,t)\in \bigl(\Omega\setminus\overline{D_j}\bigr)\times(0,T],\\[5pt]
\dfrac{\partial \rho_j^{(1)}}{\partial t}-\nu\,\Delta \rho_j^{(1)}-g_0\,\Delta u_j^{(1)}
&= 0,
&& (x,t)\in \bigl(\Omega\setminus\overline{D_j}\bigr)\times(0,T],\\[5pt]
v_j^{(1)} = 0,\quad \rho_j^{(1)}& = 0,
&& \text{on } \partial D_j\times(0,T],\\[5pt]
\partial_{\nu} v_j^{(1)} = 0,\quad \partial_{\nu} \rho_j^{(1)} &= 0,
&& \text{on } \Sigma_2,\\[5pt]
v_j^{(1)}\Big|_{t=T} = 0,\quad \rho_j^{(1)}\Big|_{t=0}&=g_1>0,
&& x\in \Omega\setminus\overline{D_j}.
\end{aligned}
\right.
\end{equation}
Combining the first equation with the homogeneous boundary data, one readily deduces $v_j^{(1)} \equiv 0$ on $\overline{Q}$, so that $\Delta v_j^{(1)} = 0$ identically.
System \eqref{mN121} thus simplifies to a single scalar heat equation for $\rho_j^{(1)}$:
\begin{equation}\label{eq:reduced-heat-m1}
\begin{cases}
\begin{aligned}
\dfrac{\partial \rho_j^{(1)}}{\partial t} - \nu\,\Delta \rho_j^{(1)} &= 0,
&& \text{in } (\Omega\setminus\overline{D_j}) \times (0,T],\\[5pt]
\rho_j^{(1)} &= 0,
&& \text{on } \partial D_j \times (0,T],\\[5pt]
\partial_\nu \rho_j^{(1)} &= 0,
&& \text{on } \Sigma_2,\\[5pt]
\rho_j^{(1)}\Big|_{t=0} &= g_1 > 0,
&& \text{in } \Omega\setminus\overline{D_j}.
\end{aligned}
\end{cases}
\end{equation}
Let $G$ stand for the connected component of $\Omega \setminus (\overline{D_1} \cup \overline{D_2})$ that contains $\partial\Omega$ as a portion of its boundary, and define the difference function $V := \rho_1^{(1)} - \rho_2^{(1)}$.
By the equality of the measurement maps stated in \eqref{cod2}, the function $V$ satisfies
\begin{equation}\label{eq:difference-V-system}
\begin{cases}
\begin{aligned}
\partial_t V - \nu\Delta V &= 0, \qquad && \text{in } G \times (0,T],\\[5pt]
V &= 0, && \text{on } S_T,\\[5pt]
\partial_\nu V &= 0, && \text{on } \Sigma_2,\\[5pt]
V\big|_{t=0}&= 0, && \text{in } G.
\end{aligned}
\end{cases}
\end{equation}
As shown in Theorem \ref{Nthm2}, the unique continuation for linear parabolic equations coupled with the maximum principle gives \(D_1=D_2\).

\emph{Step II.} We work under the condition $D_1 = D_2$, which has been established in Step I, and proceed to carry out the second-order linearization of the system.
For each index $j = 1,2$, the second-order variations $v_j^{(1,2)}$ and $\rho^{(1,2)}$ satisfy the coupled system below:
\begin{equation}\label{umN21}
\left\{
\begin{aligned}
-\dfrac{\partial v_j^{(1,2)}}{\partial t}-\nu\,\Delta v_j^{(1,2)}
+\Bigl(\nabla v^{(1)}\Bigr)\cdot \Bigl(\nabla v^{(2)}\Bigr)
&= K_j^{(2)}\, \rho^{(1)} \rho^{(2)},
&& (x,t)\in Q,\\[5pt]
\dfrac{\partial \rho^{(1,2)}}{\partial t}-\nu\,\Delta \rho^{(1,2)}
-\nabla\cdot\Bigl(\rho^{(1)}\nabla v^{(2)}\Bigr)
&-\nabla\cdot\Bigl(\rho^{(2)}\nabla v^{(1)}\Bigr)
-g_0\,\Delta v_j^{(1,2)}
= 0,
&& (x,t)\in Q,\\[5pt]
v_j^{(1,2)} &= 0,\quad \rho^{(1,2)} = 0,
&& \text{on } \Sigma_1,\\[5pt]
\partial_{\nu} v_j^{(1,2)} &= 0,\quad \partial_{\nu} \rho^{(1,2)} = 0,
&& \text{on } \Sigma_2,\\[5pt]
v_j^{(1,2)}\Big|_{t=T} &= 0,\quad \rho^{(1,2)}\Big|_{t=0}=0,
&& x\in \Omega\setminus\overline{D}.
\end{aligned}
\right.
\end{equation}
From Step I we know that $v^{(1)} \equiv 0$ and $v^{(2)} \equiv 0$ on $\overline{Q}$, which immediately gives $\nabla v^{(1)} = \nabla v^{(2)} = 0$.
Substituting these vanishing terms into \eqref{umN21} reduces it to
\begin{equation}\label{umNs2}
\left\{
\begin{aligned}
-\dfrac{\partial v_j^{(1,2)}}{\partial t}-\nu\,\Delta v_j^{(1,2)}
&= K_j^{(2)}\, \rho^{(1)} \rho^{(2)},
&& (x,t)\in Q,\\[5pt]
\dfrac{\partial \rho^{(1,2)}}{\partial t}-\nu\,\Delta \rho^{(1,2)}&-g_0\,\Delta v_j^{(1,2)}
= 0,
&& (x,t)\in Q,\\[5pt]
v_j^{(1,2)} &= 0,\quad \rho^{(1,2)} = 0,
&& \text{on } \Sigma_1,\\[5pt]
\partial_{\nu} v_j^{(1,2)} &= 0,\quad \partial_{\nu} \rho^{(1,2)} = 0,
&& \text{on } \Sigma_2,\\[5pt]
v_j^{(1,2)}\Big|_{t=T} &= 0,\quad \rho^{(1,2)}\Big|_{t=0}=0,
&& x\in \Omega\setminus\overline{D}.
\end{aligned}
\right.
\end{equation}
We now subtract the system for $j=2$ from the one for $j=1$, and introduce the difference function $p := v_1^{(1,2)} - v_2^{(1,2)}$.
This difference satisfies the boundary value problem
\begin{equation}\label{eq:difference-p-secondorder}
\begin{cases}
\begin{aligned}
-\partial_t z - \nu\Delta z &= \bigl(K_1^{(2)} - K_2^{(2)}\bigr) \rho^{(1)} \rho^{(2)}, \qquad && \text{in } Q,\\[5pt]
z &= 0, && \text{on } \Sigma_1,\\[5pt]
\partial_\nu z &= 0, && \text{on } \Sigma_2,\\[5pt]
z\big|_{t=T} &= 0, && \text{in } \Omega\setminus\overline{D}.
\end{aligned}
\end{cases}
\end{equation}
Let $\tilde{b}$ denote the solution to the adjoint problem introduced earlier.
Upon multiplying equation \eqref{eq:difference-p-secondorder} with \(\tilde{b}\), integrating across the space‑time domain $Q$, we further make use of Green’s identity as well as the equality coming from the available measurement data
$\int_{U_T} v_1^{(1,2)} n \,dxdt = \int_{U_T} v_2^{(1,2)} n \,dxdt$,
we obtain the integral identity
\begin{equation*}
\iint_{(0,T)\times(\Omega\setminus\overline{D})}
\big(K_1^{(2)}(x)-K_2^{(2)}(x)\big)\,
\rho^{(1)}(x,t)\,\rho^{(2)}(x,t)\,\tilde{b}(x,t)\,dxdt=0
\end{equation*}
In light of Lemma 3.1, we may take $\rho^{(2)}(x,t) = e^{-\mu t} \alpha(x,\rho)$, where $\alpha(\cdot,\rho)$ is an arbitrary mixed eigenfunction of $-\Delta$ on $\Omega\setminus\overline{D}$.
By the completeness of the eigenfunction system in $L^2(\Omega\setminus\overline{D})$, we deduce that
\begin{equation*}
\big(K_1^{(2)}(x)-K_2^{(2)}(x)\big)
\int_{0}^{T} e^{-\mu t}\, \rho^{(1)}(x,t)\,\tilde{b}(x,t)\,dt = 0
\end{equation*}
By the strong maximum principle, both $\tilde{b}$ and $\rho^{(1)}$ remain strictly positive throughout $Q$, which ensures the integral factor does not vanish.
Hence we conclude $K_1^{(2)}(x) = K_2^{(2)}(x)$ for all $x \in \Omega\setminus\overline{D}$.

\emph{Step III}. Under the induction hypothesis $K_1^{(i)}=K_2^{(i)}$ for $i=1,\dots,N-1$ ($N\ge3$), we apply an $N$-th order mixed linearization $\partial_{\epsilon_1}\cdots\partial_{\epsilon_N}$ at $\epsilon=0$ and subtract the equations for $j=1,2$ to obtain
 \[
\iint_{(0,T)\times(\Omega\setminus\overline{D})}
\big(K_1^{(N)}(x)-K_2^{(N)}(x)\big)\,
\rho^{(1)}\,\rho^{(2)}\cdots \rho^{(N)}\,\tilde{b}(x,t)\,dxdt=0
\]
With $\rho^{(2)}=e^{-\mu  t}\alpha(x,\rho)$, and noting that $\rho^{(1)},\rho^{(3)},\dots,\rho^{(N)}$ are positive in $\overline{Q}$ (as heat solutions with positive data) and $\tilde b>0$, it follows that $K _1^{(N)}(x)=K_2^{(N)}(x)$. Thus $K _1(x,\rho)=K _2(x,\rho)$ on $(\Omega\setminus \overline{D})\times\mathbb{R}^+$.

\end{proof}

\section*{Acknowledgments}
\addcontentsline{toc}{section}{Acknowledgments}
The work of M. Ding is supported by the National Natural Science Foundation of China under Grant No.1250012409 and the Fundamental Research Funds for the Central Universities under Grant D5000240301. The work of Hongyu Liu is supported by the Hong Kong RGC General Research Funds (Projects 11311122, 11300821, and 12301420), the NSFC/RGC Joint Research Fund (Project N\_CityU101/21), and the France‑Hong Kong ANR/RGC Joint Research Grant (A\_City203/19). The work of Guang‑Hui Zheng is supported by the National Natural Science Foundation of China (Grant No.12271151) and the Natural Science Foundation of Hunan Province (Grant No.2020JJ4166).
\bibliographystyle{abbrv}
\bibliography{PMFG}
\end{document}